\documentclass[11pt]{article}
\usepackage{amsthm, amsmath, amssymb, amsfonts, url, booktabs, tikz, setspace, fancyhdr, amsbsy}
\usepackage[margin = 0.8in]{geometry}
\usepackage{esint}
\usepackage{verbatim}
\usepackage{comment}
\usepackage{indentfirst}

\newtheorem{theorem}{Theorem}[section]

\newtheorem{lemma}[theorem]{Lemma}

\theoremstyle{definition}

\newtheorem{remark}[theorem]{Remark}

\newcommand{\norm}[1]{\left\lVert#1\right\rVert}

\newcommand{\R}{\mathbb{R}}
\newcommand{\N}{\mathbb{N}}

\newcommand{\lt}{\left}
\newcommand{\rt}{\right}

\newcommand{\uu}{\boldsymbol{u}}
\newcommand{\BB}{\boldsymbol{B}}
\newcommand{\CC}{\boldsymbol{C}}
\newcommand{\DD}{\boldsymbol{D}}
\newcommand{\EE}{\boldsymbol{E}}

\newcommand{\vv}{\boldsymbol{v}}
\newcommand{\SSS}{\boldsymbol{S}}
\newcommand{\GGG}{\boldsymbol{G}}

\newcommand{\ww}{\boldsymbol{w}}

\newcommand{\dx}{\,\mathrm{d}x}

\newcommand{\dt}{\,\mathrm{d}t}
\newcommand{\ds}{\,\mathrm{d}s}

\newcommand{\pa}{\partial}
\newcommand{\fff}{\boldsymbol{f}}
\newcommand{\ffg}{\boldsymbol{g}}
\newcommand{\ffpsi}{\boldsymbol{\psi}}

\newcommand{\bfg}{{\bf g}}

\newcommand{\bfG}{{\bf G}}

\newcommand{\bfS}{{\bf S}}

\newcommand{\calD}{\mathcal D}

\numberwithin{equation}{section}

\def\ocirc#1{\ifmmode\setbox0=\hbox{$#1$}\dimen0=\ht0
    \advance\dimen0 by1pt\rlap{\hbox to\wd0{\hss\raise\dimen0
    \hbox{\hskip.2em$\scriptscriptstyle\circ$}\hss}}#1\else
    {\accent"17 #1}\fi}

\begin{document}

\title{On the existence of strong solutions for unsteady motions of incompressible chemically reacting generalized Newtonian fluids}

\author{
Kyueon Choi\thanks{Department of Mathematics, Yonsei University, Seoul, Republic of Korea. Email: \tt{gyueon@yonsei.ac.kr}},
~Kyungkeun Kang\thanks{Department of Mathematics, Yonsei University, Seoul, Republic of Korea. Email: \tt{kkang@yonsei.ac.kr}}
~and
~Seungchan Ko\thanks{Department of Mathematics, Inha University, Incheon, Republic of Korea. Email: \tt{scko@inha.ac.kr}}
}

\date{}

\maketitle

~\vspace{-1.5cm}

\begin{abstract}
We consider a system of nonlinear partial differential equations modeling the unsteady motion of an incompressible generalized Newtonian fluid with chemical reactions. The system consists of the generalized Navier-Stokes equations with power-law type viscosity with a power-law index depending on the concentration, and the convection-diffusion equation which describes chemical concentration. This
system of partial differential equations arises in the mathematical models describing the synovial fluid which can be found in the cavities of movable joints. We prove the existence of a global strong solution for the two and three-dimensional spatially periodic domain, provided that the power-law index is greater than or equal to $(d+2)/2$ where $d$ is the dimension of the spatial domain. Moreover, we also prove that such a solution is unique under the further assumption that $p^+ < \frac{3}{2} p^-$ for the two-dimensional case and $p^+ < \frac{7}{6}p^-$ for the three-dimensional case, where $p^-$ and $p^+$ are the lower and upper bounds of the power-law index $p(\cdot)$ respectively.
\end{abstract}

\noindent{\textbf{Keywords}{: Non-Newtonian fluid, convection-diffusion equation, synovial fluid, strong solution, maximal regularity, variable power-law index}

\smallskip

\noindent{\textbf{AMS Classification:} 76D03, 35Q35, 35Q92, 76R50

\section{Introduction}

In this paper, we will prove the existence of a strong solution to a system of partial differential equations (PDEs) describing the unsteady motion of an incompressible chemically reacting flow in two and three dimensions. More precisely, we would like to know whether there exists a strong solution for the following system of PDEs:
\begin{align}
\pa_t\vv + {\rm div}\,(\vv \otimes \vv) - {\rm div}\,\SSS(c, \DD \vv)  &= - \nabla \pi  + \fff , \label{main_sys} \\
{\rm div}\, \vv &= 0, \\ 
\pa_t c + {\rm div}\,(c \vv) -\Delta c  &= -{\rm div} \, \ffg \label{concentration}
\end{align}
in $Q_T :=\Omega\times(0,T)$ where $\Omega=[0,1]^d$ with $d=2,3$. Here $\vv: Q_T \to \R^d$, $\pi: Q_T \to \R$ and $c: Q_T \to \R_+$ represent the velocity field, the scalar pressure and the concentration distribution respectively. In addition, $\fff: Q_T\to \R^d$ denotes a given external forcing term for the velocity and $\ffg$ represents the term denoting the chemical concentration source. Here $\DD \vv$ is the symmetric part of the velocity gradient $\nabla \vv$, i.e., $\DD \vv = \frac{1}{2}\lt(\nabla \vv + (\nabla \vv)^T\rt)$ and $\SSS(c,\DD \vv): Q_T \to \R^d\times \R^d$ is the Cauchy stress tensor. The given system \eqref{main_sys}-\eqref{concentration} consists of the generalized incompressible Navier--Stokes equations and the convection-diffusion equation for the concentration.
In this setting, we consider the initial-boundary value problem described by \eqref{main_sys}-\eqref{concentration}. On the boundary, we assume that all involved quantities regarding the velocity $\vv$ and the concentration $c$ are periodic with respect to the domain $[0,1]^d$. In addition, we impose the following initial conditions
\begin{equation*} 
\vv(x,0)=\vv_0(x)\quad\mbox{and}\quad c(x,0)=c_0(x)\quad\mbox{in } \Omega.
\end{equation*}
For the Cauchy stress tensor $\SSS$, we consider the case when it can be represented as the form
\begin{equation*}
\bfS(c,\DD\vv)=2\nu(c,|\DD\vv|)\DD\vv,
\end{equation*}
where $\nu$ is of the power-law-like viscosity
\begin{equation*}
    \nu(c, |\DD\vv|) = \nu_0 (1 + |\DD\vv|^2)^\frac{p(c)-2}{2}
\end{equation*}
with some positive constant $\nu_0>0$. Moreover, we assume that the exponent $p: \R \to \R$ is a Lipschitz continuous function with $1 < p^- \leq p(\cdot)\leq p^+ <\infty$. Then it is straightforward to verify that the following properties hold: there exist positive constants $K_1$, $K_2$ and $K_3$ such that for any $(c,\DD)\in\R\times\R^{d\times d}$ and $\BB \in \mathbb{R}^{d \times d}$, there holds
\begin{equation} \label{P1}
\begin{split}
\frac{\pa \SSS(c,\DD)}{\pa \DD}:(\BB\otimes\BB) &\ge K_1(1+|\DD|^2)^\frac{p(c)-2}{2}|\BB|^2, \\
\lt|\frac{\pa\SSS(c,\DD)}{\pa \DD}\rt| &\le K_2(1+|\DD|^2)^\frac{p(c)-2}{2}, \\
\lt|\frac{\pa\SSS(c,\DD)}{\pa c}\rt| &\le K_3(1+|\DD|^2)^\frac{p(c)-1}{2}\log(2+|\DD|),
\end{split}
\end{equation}
with the notation $(\BB\otimes\BB)_{ijkh}=B_{ij}B_{kh}$. Furthermore, we can also deduce from the above properties that for any $c \in \mathbb{R}$ and $\DD_1, \DD_2 \in \mathbb{R}^{d \times d}$, there exists a positive constant $K_4$ such that
\begin{equation} \label{P2}
    \left(\SSS(c,\DD_1)-\SSS(c,\DD_2)\right):(\DD_1 - \DD_2) \ge K_4 (1+ |\DD_1|^2+ |\DD_2|^2)^\frac{p(c)-2}{2} |\DD_1-\DD_2|^2.
\end{equation}

This model generalizes the power-law-like non-Newtonian fluids, which was first investigated in the late 1960s by Ladyzhenskaya and Lions independently in \cite{Lady1967}, \cite{Lady1969} and \cite{Lion1969}. They constructed weak solutions of the steady-state model for $p \ge \frac{3d}{d+2}$ as well as the unsteady problem for $p \ge \frac{3d+2}{d+2}$ by using the monotone operator theory. Further in \cite{Malek2003} and \cite{Diening2010}, the authors relaxed the condition of $p$ to $p > \frac{2d}{d+2}$ for both stationary and nonstationary cases by applying the so-called Lipschitz truncation technique. 

However, laboratory experiments have revealed that in several non-Newtonian fluid flow models, the power-law index is not a fixed constant
and often varies. More precisely, in some cases, it was found that the power-law exponent can be represented as $p(\cdot) = p(\EE(x))$, where $\EE$ is an electric field described by the quasi-static Maxwell equations, which is called the electrorheological fluids \cite{Ruzicka1996, Ruzicka2004}. Furthermore, \cite{Ruzicka2000} studied the existence theory for electrorheological fluids. In \cite{Diening2008}, the authors improved the condition for $p(\cdot)$ by using the Lipschitz truncation method. To be more specific, the paper shows the existence of weak solutions for sufficiently regular $p(\cdot)$ with $1 < p^- \le p(\cdot) \le p^+ < \infty$ and $p^- > \frac{2d}{d+2}$ for the stationary model. More recently, the model consisting of the generalized Newtonian fluids with concentration-dependent viscosity coupled with the convection-diffusion equation was studied in \cite{Bulicek2009}. The authors considered the concentration effect as a scaling factor of the viscosity, i.e., $\nu(c,\DD) \sim f(c)\tilde{\nu}(\DD)$, and the power-law exponent remained constant. However, the experimental results for the biological fluids, e.g. synovial fluid, showed that the concentration affects the level of shear-thinning effects. For a detailed explanation of the mathematical modeling of that fluid, we refer to \cite{Hron2010} and \cite{Pustejovska2012}. These fluids are called the chemically reacting fluids (CRF) and the existence theory of weak solutions for stationary case was developed in \cite{Bulicek2013} when $p^- > \max\{\frac{3d}{d+2}, \frac{d}{2}\}$ and improved to $p^- > \frac{d}{2}$ in \cite{Bulicek2014}. Furthermore, in \cite{Ko2022}, the author proved the existence of weak solutions when $p^- > \frac{d+2}{2}$ for the non-stationary case. From the viewpoint of computational mathematics, \cite{ko2d,ko3d} constructed the finite element approximation of the weak solutions and performed the convergence analysis of the numerical methods. On the other hand, for strong solutions, there are only a few results known. The global existence of strong solutions for $p\ge\frac{3d+2}{d+2}$ was shown in \cite{Bellout1994} and \cite{Malek1993}. Furthermore, the local existence of strong solutions for large data and the global existence of strong solutions for small data can be found in \cite{Malek1995} when $p>\frac{3d-4}{d}$. In \cite{Diening2005}, the authors established the local existence of strong solutions for large data for $\frac{7}{5} < p \le 2$ in three dimensions, which can be extended to the case of electrorheological fluids concerning the variable power-law index. Moreover, in \cite{kang_1, kang_2}, the authors proved the existence of regular solutions for the generalized Newtonian fluids under certain conditions, and \cite{kang_3, kang_4} considered the non-Newtonian fluid flow models coupled with other equations and established the existence of solutions with higher regularity. For the model of CRF, the existence of classical solutions was shown for $d=2$ in \cite{Bulicek2019}.

To the best of our knowledge, there are no results about the existence of strong solutions for the unsteady CRF model. A significant difference compared to the electrorheological fluids is the appearance of $\nabla c$ and $\partial_t c$ when we differentiate the Cauchy stress tensor in the Galerkin method, which requires us to show the higher integrability of the derivatives of the concentration. Since the convection-diffusion equation can be regarded as the heat equation with a source term, we apply the maximal regularity theory to obtain the desired results. The rest of this paper is organized as follows. We first introduce the auxiliary tools for the proofs of our main theorems in Section 2. In the following section, we construct the Galerkin approximation as a baseline system and prove the higher regularity of $\nabla c$. In Section 4 and Section 5, we derive some uniform estimates including the $W^{2,2}$-estimate for the velocity $\vv$ and $L^2$-estimate of $\pa_t \vv$ both for the two-dimensional case and the three-dimensional case. Then in Section 6, we prove the existence of a strong solution through the limiting process based on the uniform estimates obtained in the previous sections. Further, we prove the uniqueness of the strong solution under the additional assumption in Section 7 and we will discuss the concluding remark in Section 8.

\section{Preliminaries and Main Result}

In this section, we shall introduce our main result along with some preliminaries which will be used throughout the paper. For two vectors $\boldsymbol{a}$, $\boldsymbol{b}$, $\boldsymbol{a}\cdot \boldsymbol{b}$ means the dot product and for two tensors $\mathbb{A}$, $\mathbb{B}$, $\mathbb{A}:\mathbb{B}$ denotes the scalar product. In addition, $C$ signifies a generic positive constant, which may vary at each appearance. For $k\in\mathbb{N}\cup\{0\}$ and $1\leq q\leq\infty$, we denote the standard Lebesgue and Sobolev spaces by $L^q(\Omega)$ and $W^{k,q}(\Omega)$ respectively. Moreover, for the sake of simplicity, we shall write $\|\cdot\|_p=\|\cdot\|_{L^p(\Omega)}$ and $\|\cdot\|_{k,p}=\|\cdot\|_{W^{k,p}(\Omega)}$. For a Banach space $X$, we define the Bochner space $L^q(0,T;X)$ by the family of functions with the finite Bochner-type norm
\begin{align*}
    & \|u\|_{L^q(0,T;X)} := \begin{cases}
        \bigg(\displaystyle \int_0^T \|u(t)\|_X^q \dt \bigg)^\frac{1}{q} & \text{if } 1 \le q < \infty, \\
        \displaystyle\sup_{t\in(0,T)} \|u(t)\|_X & \text{if } q = \infty,
    \end{cases}
\end{align*}
with an abbreviation $L^q(Q_T)=L^q(0,T;L^q(\Omega))$. Since we will deal with the spatially periodic functions, we shall restrict ourselves to the case of the zero mean-value functions over $\Omega$. In this perspective, we will use the following spaces which are frequently used in the study of incompressible fluid flow problems: For $k \in \mathbb{N} \cup \{0\}$ and $1 \le q < \infty$,
\begin{align*}
&\mathcal{V}:= \{ \varphi \in C^\infty(\Omega) : \langle \varphi,1\rangle = 0,\,\,{\rm div}\, \varphi =0 \}, \\
&W^{k,q}_{\rm div} (\Omega) := \overline{\mathcal{V}}^{\|\cdot\|_{k,q}}.
\end{align*}

Next, let us introduce the function space with variable exponent. Let $\mathcal{P}$ be a set of all measurable functions $p: Q_T \to [1,\infty]$ and we call the function $p \in \mathcal{P}(Q_T)$ a variable exponent. Furthermore, we define $p^- := \inf_{Q_T} p(x,t)$ and $p^+ := \sup_{Q_T} p(x,t)$.  Then the variable-exponent Lebesgue space is defined as
\[
    L^{p(\cdot)}(Q_T) := \bigg\{ u \in L^1_{\rm loc} (Q_T): \int_{Q_T} |u(x,t)|^{p(x,t)} \dx \dt < \infty \bigg\} 
\]
equipped the corresponding Luxembourg norm
\[
\|u\|_{L^{p(\cdot)}(Q_T)} := \inf \bigg\{ \lambda>0 : \int_{Q_T} \bigg| \frac{u(x,t)}{\lambda} \bigg|^{p(x,t)} \dx \dt \le 1 \bigg \}.
\]
If $1<p^- \le p^+ <\infty$, this space becomes a reflexive and separable Banach space. Furthermore, H\"older's inequality and Young's inequality hold with the variable exponents $1= \frac{1}{p(\cdot)} + \frac{1}{p'(\cdot)}$. These properties play an important role as tools for the analysis of the CRF under consideration. See, e.g.,  \cite{diening} for more information on the variable-exponent Lebesgue and Sobolev spaces. 

Next, We shall present some auxiliary results which will be needed in the later analysis. We begin with the following well-known lemma, called Korn's inequality. We refer to p.196 in \cite{malek} for the details.
\begin{lemma}[Korn's inequality] \label{KI}
For $1<q<\infty$, assume that $\uu \in W^{1,q}(\Omega)$. Then there holds
\[
    \|\nabla \uu\|_q \le C(q,\Omega) \|\DD \uu\|_q,
\]
where $\Omega$ is $d$-dimensional torus.
\end{lemma}
We also introduce the Gagliardo--Nirenberg interpolation theorem. See \cite{niren} for the details.

\begin{lemma}[Gagliardo--Nirenberg inequality on bounded domain]\label{GNI}
Let $\Omega\subset\R^d$ be a measurable, bounded, open and connected domain with a Lipschitz boundary. Let $1\le q,r\le \infty$, $j$ and $m$ be non-negative integers, and $p \ge 1$ such that the relation
\[
\frac{1}{p}= \frac{j}{d} + \theta\lt(\frac{1}{r}-\frac{m}{d}\rt)+\frac{1-\theta}{q}, \quad \frac{j}{m} \le \theta \le 1
\]
holds. Then we have
\[
\|\calD^j u\|_{p} \le C\|\calD^m u\|_{r}^{\theta}\|u\|_{q}^{1-\theta} + C\|u\|_{\sigma},
\]
where $u \in L^q(\Omega)$ with $\calD^m u\in L^r(\Omega)$ and $\sigma>0$ is arbitrary.

\end{lemma}

Since we will use the Galerkin method, we need to show the solvability of the system of ordinary differential equations (ODEs) for approximate solutions. The following lemma ensures the existence of $\CC : (t_0-\delta,t_0+\delta) \to \mathbb{R}^d$ solving the following system of ODEs (see, e.g., \cite{Zeidler1990}):
\begin{equation} \label{systemODE}
\begin{split}
     \frac{\rm d}{\dt} \CC(t) &= \GGG(t,\CC(t)), \quad \forall t \in (t_0-\delta,t_0 + \delta), \\
     \CC(t_0) &= \CC_0,
\end{split}
\end{equation}
where $\CC_0 \in \mathbb{R}^d$. 
\begin{lemma}[Carath\'eodory's theorem]\label{aux:cara}
Suppose that $\GGG : (t_0-\delta,t_0+\delta) \times B_{\epsilon}(\CC_0) \to \R^d$ satisfies the following properties: if we write $\bfG=\{G_i\}_{i=1}^d$,
\begin{itemize}
\item $G_i(\cdot, \CC)$ is measurable for all $i=1,\ldots,d$ and for all $\CC \in B_{\epsilon}(\CC_0)$.
\item $G_i(t,\cdot)$ is continuous for a.e. $t\in(t_0-\delta,t_0+\delta)$.
\item $|\GGG(t,\CC)|\le F(t)$ for all $(t,\CC)\in(t_0-\delta,t_0+\delta)\times B_{\epsilon}(\CC_0)$ and $F(t)$ is integrable on $(t_0-\delta,t_0+\delta)$.
\end{itemize}
Then there exist a $\delta' \in (0,\delta)$ and a continuous function $\CC: (t_0 - \delta',t_0 + \delta') \to \mathbb{R}^d$ such that $\frac{\rm d}{\dt} \CC(t)$ exists for a.e. $t \in (t_0 - \delta',t_0 + \delta')$ and $\CC$ solves \eqref{systemODE}.
\end{lemma}

The following lemma is so-called Gr\"onwall's inequality which is useful for parabolic problems. The details can be found in the Appendix of \cite{evans}.

\begin{lemma}[Gr\"onwall's inequality]
    Let $\eta$ be a nonnegative and absolutely continuous function on $[0,T]$ which satisfies for a.e. $t$ that
    \[
    \eta'(t) \le \phi(t) \eta(t) + \psi(t),
    \]
    where $\phi(t)$ and $\psi(t)$ are nonnegative, integrable functions on $[0,T]$. Then we have
    \[
    \eta(t) \le \exp\left(\int_0^T \phi(s) \ds\right) \left[ \eta(0) + \int_0^t \psi(s) \ds \right]
    \]
    for all $0 \le t \le T$.
\end{lemma}

We will also use the following well-known compactness theorem in Bochner spaces (see, e.g., \cite{Lion1969}).

\begin{lemma}[Aubin--Lions Lemma]\label{ALL}
Let $V_1$, $V_2$ and $V_3$ be reflexive and separable Banach spaces such that 
\begin{equation*}
V_1\hookrightarrow\hookrightarrow V_2 \quad\mbox{and}\quad V_2 \hookrightarrow V_3.
\end{equation*}
Then for $1<p,q<\infty$, the space $\{v \in L^p(0,T;V_1) : \pa_t \vv \in L^{q}(0,T;V_3)\}$ is compactly embedded into $L^p(0,T;V_2)$.

\end{lemma}
Furthermore, as mentioned earlier, we need the following lemma, known as the maximal-regularity estimate, to obtain the higher integrability of the derivatives of the concentration. For the details, see for example, Theorem 5.4 and Theorem D.12 of \cite{james} and p.288-291 in \cite{lady}.

\begin{lemma}[Maximal regularity] \label{MR}
Let $\Omega$ be a $d$-dimensional torus. If $u$ is the solution of 
\[
\begin{cases}
     \pa_t u - \Delta u = h, \\
     u(x,0) = u_0(x),
\end{cases}
\]
with $h \in L^r(0,T;L^\ell(\Omega))$ and $u_0 \in W^{2,\ell}(\Omega)$ for $1<r,\,\ell<\infty$. Then there holds
\[
\|\pa_t u\|_{L^r(0,T;L^\ell(\Omega))} + \|\Delta u \|_{L^r(0,T;L^\ell(\Omega))} \le C \|h\|_{L^r(0,T;L^\ell(\Omega))} + C\|\Delta u_0\|_{L^\ell(\Omega)}.
\]
\end{lemma}

\begin{remark}
    In the reference \cite{james}, the authors considered the case of $u_0 = 0$. We have extended this result to the case of non-homogeneous initial condition in a straightforward manner, by using the linearity of the heat equation and the lifting argument. Therefore the condition for $u_0$ presented in the above lemma might not be optimal, and one may be able to improve the sufficient condition of $u_0$ in Lemma \ref{MR}.
\end{remark}

Finally, we shall make a comment on the well-known Calder\'on--Zygmund inequality on the torus (see, e.g., Theorem B.7 in \cite{james}), which will be used to prove the uniqueness result later.
\begin{lemma}[Calder\'on--Zygmund inequality] \label{cz}
    Let $\Omega$ be a $d$-dimensional torus. Assume that $-\Delta u = f$ in $\Omega$ and $f \in L^q(\Omega)$ for some $1<q<\infty$. Then there holds
    \[
    \|\pa_j \pa_k u\|_q \le C(q) \|f\|_q.
    \]
\end{lemma}

Now we are ready to state our main theorem, concerning the existence of a strong solution of the system under consideration.

\begin{theorem} \label{maintheorem}
Let $\Omega = [0,1]^d$ with $d=2,3$. Suppose that $p: \R \to \R$ is a Lipschitz continuous function with $\frac{d+2}{2} \le p^- \leq p(\cdot)\leq p^+ <\infty$. Assume further that $\vv_0 \in W_{\rm div}^{1,p^+}(\Omega)$, $\fff \in L^2(0,T;L^2(\Omega))$, $c_0 \in W^{2,q}(\Omega)$ and $\ffg \in L^q(0,T;W^{1,q}(\Omega))$ with $q>2d$. Then, there exists a strong solution $(\vv,\pi, c)$ such that 
\begin{align*}
    & \vv \in C([0,T];L^2(\Omega)) \cap L^\infty(0,T;W_{\rm div}^{1,2}(\Omega)) \cap L^2(0,T; W^{2,2}(\Omega)),\\
    & |\nabla \vv|^\frac{p(c)}{2} \in L^2(0,T;W^{1,2}(\Omega)) \cap L^\infty(0,T;L^2(\Omega)),\\
    &\partial_t \vv \in L^2(0,T;L^2(\Omega)), \quad \nabla\pi\in L^{\frac{6}{5}}(0,T;L^{\frac{6}{5}}(\Omega)),\\
    & c \in C([0,T];L^2(\Omega)) \cap L^\infty(0,T;W^{1,q}(\Omega)) \cap L^2(0,T;W^{2,2}(\Omega)),\\
    & |\nabla c |^\frac{q}{2} \in L^2(0,T;W^{1,2}(\Omega)),\quad  \pa_t c \in L^2(0,T;L^2(\Omega)),
\end{align*}
which satisfy the given equations \eqref{main_sys}-\eqref{concentration} for a.e. in $Q_T$.

\end{theorem}

As our second result, we will also prove the following theorem regarding the uniqueness of the strong solution in the above theorem.
\begin{theorem} \label{uniqueness theorem}
    Assume that $p^+ < \frac{3}{2} p^-$ for the two-dimensional case and $p^+<\frac{7}{6}p^-$ for the three-dimensional case. Then the strong solution constructed in Theorem \ref{maintheorem} is unique. 
\end{theorem}

\section{Galerkin approximation}

\subsection{Galerkin approximation} \label{best}
Let us first construct the Galerkin approximation for the given system \eqref{main_sys}-\eqref{concentration}.  Let $\{\ww_i\}_{i=1}^{\infty}$ denotes the set of eigenvectors of the Stokes operator with $\langle \ww_i,1\rangle=0$ for all $i\in\mathbb{N}$. Thanks to the periodic boundary condition, we can formally write
\begin{equation} \label{vbasis}
-\Delta \ww_i = \lambda_i \ww_i \quad \mbox{in }\Omega,
\end{equation} 
where $\lambda_i$ is the eigenvalue of the Stokes operator corresponding $\ww_i$. We further note that $\{\ww_i\}_{i=1}^\infty$ forms a basis of $W^{s,2}_{\rm div}(\Omega)$ for $s \in \mathbb{N}$ (see Theorem 4.11, p.290 in \cite{malek}). Moreover, we choose $s>1+\frac{d}{2}$ so that $W^{s,2}_{\rm div}(\Omega) \hookrightarrow W^{1,\infty}_{\rm div}(\Omega)$. For the concentration equation, we shall consider $\{z_i\}_{i=1}^{\infty}$ forming a basis of $W^{1,2}(\Omega)$ such that  $\int_{\Omega} z_i z_j \dx = \delta_{ij}$. Then for each $n$, $m\in\N$, we seek for a pair $(\vv^{n,m}, c^{n,m})$ given by
\[
\vv^{n,m} := \sum_{i=1}^n \alpha_i^{n,m}(t)\ww_i, \quad c^{n,m} := \sum_{j=1}^m \beta_j^{n,m}(t)z_j
\]
satisfying the following system of ODEs: 

\begin{equation}\label{aux1}
\begin{split}
\int_{\Omega} \pa_t\vv^{n,m} \cdot \ww_i \dx - \int_{\Omega} (\vv^{n,m}\otimes\vv^{n,m}):\nabla \ww_i \dx + \int_{\Omega} \SSS^{n,m}: \DD \ww_i \dx = \int_{\Omega} \fff \cdot \ww_i \dx,
\end{split} 
\end{equation}

\begin{equation}
\vv^{n,m}(\cdot,0) = P^n\, \vv_0,
\end{equation}

\begin{equation}
\int_{\Omega} \pa_t c^{n,m} z_j \dx - \int_{\Omega}  c^{n,m} \vv^{n,m}\ \cdot \nabla z_j \dx + \int_{\Omega} \nabla c^{n,m} \cdot \nabla z_j \dx = \int_{\Omega} \ffg \cdot \nabla z_j \dx,
\end{equation}

\begin{equation}\label{aux4}
c^{n,m}(\cdot,0)=P^m\, c_0
\end{equation}
for all $i=1,\ldots,n$ and $j=1,\ldots,m$, where $\SSS^{n,m}:= \SSS(c^{n,m},\DD \vv^{n,m})$, and $P^n$, $P^m$ denote the $L^2$-orthogonal projections onto $A_n:={\rm span}\{\ww_1,\ldots,\ww_n\}$ and $B_m:={\rm span}\{z_1,\ldots,z_m\}$ respectively. Based on Lemma \ref{aux:cara}, it is straightforward to show the existence of solutions to the discrete equations \eqref{aux1}-\eqref{aux4} at least for a short time interval, and the solution can be extended to the whole time interval $(0,T)$ using the typical uniform estimates. Thanks to the argument used in \cite{Ko2022}, we can take the limit $m \to \infty$ to get solutions $\vv^n=
\sum_{i=1}^n\alpha_i^n(t)\ww_i\in A_n$ and $c^n\in L^2(0,T;W^{1,2}(\Omega))$ with $\pa_t c^n \in L^2(0,T;W^{-1,2}(\Omega))$ satisfying the following intermediate problem: for all $i=1,\ldots,n$, there holds
\begin{equation}\label{aux11}
\begin{split}
\int_{\Omega} \pa_t\vv^n \cdot \ww_i \dx - \int_{\Omega} (\vv^n\otimes\vv^n):\nabla \ww_i \dx + \int_{\Omega} \SSS^n: \DD \ww_i \dx = \int_{\Omega} \fff \cdot \ww_i \dx,
\end{split} 
\end{equation}

\begin{equation*}
\vv^n(\cdot,0) = P^n \, \vv_0,
\end{equation*}

\begin{equation}\label{aux33}
\int_{\Omega} \pa_t c^n \varphi \dx - \int_{\Omega}  c^n \vv^n\ \cdot \nabla \varphi \dx + \int_{\Omega} \nabla c^n\cdot \nabla \varphi \dx = \int_{\Omega} \bfg \cdot \nabla \varphi \dx,
\end{equation}

\begin{equation*}
c^n(\cdot,0)= c_0
\end{equation*}
for arbitrary $\varphi \in W^{1,2}(\Omega)$ with
$\SSS^n:= \SSS(c^n,\DD \vv^n)$.
Multiplying the $i$-th equation of \eqref{aux11} by $\alpha_i^n$ and taking the sum over $i=1,\ldots,n$, we obtain
\[
\frac{1}{2}\frac{\rm d}{{\rm d}t}\int_{\Omega} |\vv^n|^2\dx + \int_{\Omega} \SSS^n:\DD \vv^n \dx = \int_{\Omega} \fff\cdot \vv^n\dx,
\]
which leads us to derive the estimate
\begin{equation} \label{best1}
    \sup_{0\le t \le T} \norm{\vv^n}_2^2 + \int_0^T \int_{\Omega} \left(|\nabla \vv^n|^{p(c^n)} + |\SSS^n|^{(p(c^n))'} \right)\dx\dt \le C,
\end{equation}
where the constant on the right-hand side is independent of $n\in\mathbb{N}$.
If we take $\varphi = c^n$ in \eqref{aux33}, we have
\[
\frac{1}{2}\frac{\rm d}{{\rm d}t}\int_{\Omega} |c^n|^2 \dx + \int_{\Omega} |\nabla c^n|^2 \dx = \int_{\Omega} \ffg \cdot \nabla c^n \dx,
\]
from which we obtain the following uniform estimate
\begin{equation} \label{best2}
\sup_{0 \le t \le T} \norm{c^n}_2^2 + \int_0^T \int_{\Omega} |\nabla c^n|^2 \dx \dt \le C.
\end{equation}
Before proceeding further for $n \to \infty$, we need more uniform estimates regarding the higher regularity of $\nabla c^n$ which will be subsequently investigated in the following sections. Henceforth, we will omit the index $n$ for the sake of convenience.

\subsection{Higher integrability of $\nabla c$} \label{hgradc}

To begin with, we shall prove the higher integrability of $\nabla c$. Indeed, we aim to obtain
\begin{equation}\label{high_int_c}
\nabla c \in L^\infty(0,T;L^q(\Omega))
\end{equation} 
for $q>2d$. To do this, since $C^\infty(\Omega)$ is dense in $W^{1,2}(\Omega)$, we can formally take $\varphi = - \nabla \cdot ( |\nabla c|^{q-2} \nabla c)$ in \eqref{aux33} and consequently it follows that
\begin{align*}
\frac{1}{q}\frac{\rm d}{{\rm d}t}\int_{\Omega} |\nabla c|^q \dx + \int_{\Omega} |\nabla c|^{q-2}|\nabla^2 c|^2 \dx \le \int_{\Omega} |\nabla \vv||\nabla c|^q \dx + \int_\Omega {\rm div}\, \ffg \,\nabla \cdot (|\nabla c|^{q-2} \nabla c) \dx =: {\rm D_1} + {\rm D_2}.
\end{align*}
Note that ${\rm D_1}$ can be estimated as
\begin{align*}
{\rm D}_1 \le \|\nabla \vv\|_{p^-}\| |\nabla c|^{\frac{q}{2}}\|^2_{2(p^-)'} 
&\le C\|\nabla \vv \|_{p^-} \| |\nabla c |^\frac{q}{2}\|_2^{2\theta} \| \nabla (|\nabla c |^\frac{q}{2}) \|^{2(1-\theta)}_2 \\
&\le C \|\nabla \vv\|_{p^-}^\frac{1}{\theta} \| |\nabla c |^\frac{q}{2}\|_2^2 + \epsilon \| \nabla (|\nabla c |^\frac{q}{2}) \|_2^2,
\end{align*}
where we have used Young's inequality and the Gagliardo--Nirenberg interpolation inequality (Lemma \ref{GNI})
with $\frac{1}{2{(p^-)}'} = \frac{\theta}{2} + (\frac{1}{2} - \frac{1}{d})(1- \theta)$, meaning that $\theta = \frac{d-d{(p^-)}'+2{(p^-)}'}{2 {(p^-)}'}$.
Next, we estimate ${\rm D_2}$ by H\"older's inequality and Young's inequality as
\[
    {\rm D_2} \le C(q) \|\ffg\|_{1,q} \|\nabla c\|_q^\frac{q-2}{2} \bigg(  \int_\Omega |\nabla c|^{q-2}|\nabla^2 c|^2 \dx \bigg)^\frac{1}{2} \le C \|\ffg\|_{1,q}^q + C\|\nabla c\|_q^q + \epsilon \int_\Omega |\nabla c|^{q-2}|\nabla^2 c|^2 \dx.
\]
Further, note that there holds
\[
\left| \nabla (|\nabla c|^\frac{q}{2}) \right|^2 \le C(q) |\nabla c|^{q-2} |\nabla^2 c|^2.
\]
Now by combining the above inequalities, we obtain
\[
\frac{\rm d}{{\rm d}t} \int_\Omega |\nabla c |^q \dx + C\int_\Omega \left|\nabla (|\nabla c |^\frac{q}{2}) \right|^2 \dx \le C (1+\|\nabla \vv \|^{\frac{1}{\theta}}_{p^-} ) \int_\Omega |\nabla c |^q \dx + C\|\ffg\|_{1,q}^q.
\]
Recalling that $\nabla \vv \in L^{p^-}(Q_T)$ from \eqref{best1} and $\ffg \in L^q(0,T;W^{1,q}(\Omega))$ from the assumption, we need the condition $\frac{1}{\theta} \le p^-$ to apply Gr\"onwall's inequality, which is equivalent to $\frac{d+2}{2}\le p^-$. Therefore, if we assume the condition $\frac{d+2}{2}\le p^-$, by Gr\"onwall's inequality, we can conclude that
\[
\sup_{0\le t\le T}\|\nabla c(t)\|_q + \int_0^T\int_\Omega \left|\nabla (|\nabla c |^\frac{q}{2}) \right|^2 \dx\dt \le C, 
\]
provided that $\nabla c_0 \in L^{q}(\Omega)$ which comes from the assumption. Note that the above argument holds both for the two and three-dimensional domains. In the following sections, we will derive some uniform estimates for dimensions two and three separately.

\section{Uniform estimates in $2$D}

\subsection{$W^{2,2}$-estimate of $\vv$}

Now we aim to get the uniform estimate of $\|\vv\|_{L^2(0,T;W^{2,2}(\Omega))}$. If we multiply the $i$-th equation in \eqref{aux11} by $\lambda_i \alpha_i^n(t)$ and take the sum over $i=1,\ldots,n$, then in virtue of \eqref{vbasis} we get
\begin{equation} \label{W22est1}
\int_{\Omega} \pa_t \vv \cdot (-\Delta \vv) \dx - \int_{\Omega} (\vv\otimes\vv) : \nabla (-\Delta \vv) \dx + \int_{\Omega} \SSS(c,\DD \vv):\DD (-\Delta \vv) \dx = \int_{\Omega} \fff \cdot (-\Delta \vv) \dx.
\end{equation}
Note that in the two-dimensional domain, we can verify by the integration by parts that
\[
    \int_{\Omega}(\vv\otimes\vv):\nabla\Delta\vv\dx=\sum^2_{i,j,k=1}\int_{\Omega}\partial_{x_k}v_j \partial_{x_j}v_i \partial_{x_k}v_i\dx=0.
\]
Therefore, \eqref{W22est1} can be written as
\begin{equation} \label{W22est2}
\frac{1}{2}\frac{\rm d}{{\rm d}t}\int_{\Omega} |\nabla \vv|^2 \dx + {\rm I} = \int_{\Omega} \fff \cdot (-\Delta \vv) \dx,
\end{equation}
with
\[
{\rm I} := \int_{\Omega}  \frac{\pa \SSS(c,\DD)}{\pa \DD}:(\pa_{x_k}\DD\vv\otimes\pa_{x_k}\DD\vv) \dx+\int_{\Omega}  \pa_{x_k} c  \left(\frac{\pa\SSS(c,\DD)}{\pa c} : \pa_{x_k}\DD \vv\right)
  \dx 
=: {\rm E_1} + {\rm E_2}.
\]
First, due to \eqref{P1}, it follows that 
\begin{equation}\label{3_3_3}
    {\rm E_1} \ge C \int_\Omega ( 1+ |\DD \vv|^2)^\frac{p(c)-2}{2} |\nabla \DD \vv|^2 \dx.
\end{equation}
For the estimate of ${\rm E_2}$, by \eqref{P1} and H\"older's inequality, we have
\begin{equation}\label{3_3_1}
\begin{split}
    |{\rm E_2}| &\le C \int_\Omega |\nabla c | ( 1+ |\DD \vv|^2)^\frac{p(c)-1}{2} \log (2 + |\DD \vv|) |\nabla \DD \vv| \dx \\
    &= C \int_\Omega |\nabla c| \log(2+|\DD \vv|)  ( 1+ |\DD \vv|^2)^\frac{p(c)}{4} \big|( 1+ |\DD \vv|^2)^\frac{p(c)-2}{4} \nabla \DD \vv  \big| \dx \\
    &\le C \|\nabla c\|_q \| \log(2+ |\DD \vv|) \|_\alpha \| (1+ |\DD \vv|^2)^\frac{p(c)}{4} \|_4 \|  ( 1+ |\DD \vv|^2)^\frac{p(c)-2}{4} \nabla \DD \vv \|_2,
\end{split}    
\end{equation}
where 
\begin{equation} \label{alphaq2d}
    1 = \frac{1}{q} + \frac{1}{\alpha} + \frac{1}{4} + \frac{1}{2}
\end{equation}
with $q>4$ as given in the assumption of Theorem \ref{maintheorem} and $\alpha>4$ large enough. Note from \eqref{high_int_c} that we have $\sup_{0\leq t\leq T}\|\nabla c (t)\|_q \le C$. To estimate the remaining terms, let us first define 
\begin{equation} \label{eta}
    \eta := (\overline{\DD} \vv)^\frac{p(c)}{2}\quad \text{with} \quad \overline{\DD} \vv := (1 + |\DD \vv|^2)^\frac{1}{2}.
\end{equation} 
Then we can rewrite \eqref{3_3_1} as

\begin{equation} \label{E2est}
    |{\rm E_2}| \le C\| \log(2+ |\DD \vv|) \|_\alpha \| \eta \|_4 \|  (\overline{\DD}\vv)^\frac{p(c)-2}{2} \nabla \DD \vv\|_2.
\end{equation}
In order to estimate $\|\eta\|_{\beta}$, since $p \in W^{1,\infty}(\mathbb{R})$, 
\[
    |\nabla \eta| \le C |( \overline{\DD} \vv)^\frac{p(c)-2}{2} \nabla \DD \vv| + C |\nabla c (\overline{\DD} \vv)^\frac{p(c)}{2} \log (2 + |\DD \vv|)|.
\]
By H\"older's inequality and \eqref{high_int_c}, we have
\begin{equation} \label{2Dgradeta}
    \| \nabla \eta\|_2 \le C\left(\|(\overline{\DD} \vv )^\frac{p(c)-2}{2} \nabla \DD \vv \|_2 + \| \eta \|_4\|\log(2+ |\DD \vv|)\|_\alpha\right)\hspace{-0.1cm}.
\end{equation}
By the Gagliardo--Nirenberg interpolation inequality, there holds that
\begin{equation} \label{2DetaGNI}
    \|\eta\|_4 \le C \| \eta \|_2^\frac{1}{2} \|\nabla \eta\|_2^\frac{1}{2}.
\end{equation}
Therefore, we have
\[
    \|\eta\|_4 
    \le C \left(\|(\overline{\DD}\vv)^\frac{p(c)-2}{2} \nabla \DD \vv \|_2^\frac{1}{2} \|\eta\|_2^\frac{1}{2} + \|\eta\|_4^\frac{1}{2} \|\log(2+ |\DD \vv|)\|_\alpha^\frac{1}{2} \|\eta\|_2^\frac{1}{2}\right)
\]
and together with Young's inequality, we obtain
\begin{equation}\label{etabetaest}
    \|\eta\|_4 \le C \left( \|(\overline{\DD}\vv)^\frac{p(c)-2}{2} \nabla \DD \vv\|_2^\frac{1}{2} \|\eta\|_2^\frac{1}{2} + \|\log(2+|\DD \vv|)\|_\alpha \|\eta\|_2\right) \hspace{-0.1cm}.
\end{equation}
Now inserting \eqref{etabetaest} into \eqref{E2est} and using Young's inequality yield
\begin{align*}
    |{\rm E_2}| 
    &\le C\| \log(2+ |\DD \vv|)\|_\alpha \|\eta\|_2^\frac{1}{2} \| (\overline{\DD}\vv)^\frac{p(c)-2}{2} \nabla \DD \vv\|_2^\frac{3}{2} +  C\| \log(2+ |\DD \vv|)\|_\alpha^2 \|\eta\|_2 \| (\overline{\DD}\vv)^\frac{p(c)-2}{2} \nabla \DD \vv\|_2\\
    & \le \epsilon \|(\overline{\DD}\vv)^\frac{p(c)-2}{2} \nabla \DD \vv\|_2^2 + C \|\log(2+ |\DD \vv|)\|_\alpha^4 \|\eta\|^2_2.
\end{align*}
Since $\log(2+|\DD \vv|) \le C ( 1 + |\DD \vv|)^\frac{2}{\alpha}$, we know that $\|\log(2+|\DD \vv|)\|_\alpha^4 \le C + C \|\DD \vv\|_2^\frac{8}{\alpha}$.
Furthermore, as $4 < \alpha$, it is straightforward to verify
\[
    \|\log(2+|\DD \vv|)\|_\alpha^4 \le C + C \|\DD \vv\|_2^2,
\]
which leads us to the inequality
\begin{equation}\label{3_3_2}
    |{\rm E_2}| \le \epsilon \|(\overline{\DD}\vv)^\frac{p(c)-2}{2} \nabla \DD \vv\|_2^2 +C\|\eta\|^2_2+ C\|\DD\vv\|^2_2 \|\eta\|^2_2.
\end{equation}
Collecting the terms \eqref{3_3_3} and \eqref{3_3_2}, we have
\[
    {\rm I} \ge C \int_\Omega (1+ |\DD \vv|^2)^\frac{p(c)-2}{2} |\nabla \DD \vv|^2 \dx - C \|\eta\|_2^2 - C\|\DD \vv\|_2^2\|\eta\|_2^2.
\]
Therefore, from \eqref{W22est2} together with Young's inequality, we have 
\[
    \frac{1}{2} \frac{\rm d}{\dt} \|\nabla \vv\|^2_2 + C \int_\Omega (1+ |\DD \vv|^2)^{\frac{p(c)-2}{2}} |\nabla \DD \vv|^2 \dx \le  C\|\fff\|_2^2 + \epsilon \|\Delta \vv\|_2^2 + C \|\eta\|_2^2 + C\|\nabla \vv\|_2^2\|\eta\|_2^2.
\]
Note that since $p^- \ge \frac{d+2}{2} = 2$ , we know that
\[
\|\Delta\vv\|^2_2\leq C\int_\Omega |\nabla \DD\vv| ^2 \dx \le  \int_\Omega(1+ |\DD\vv|^2)^\frac{p(c)-2}{2} |\nabla \DD\vv|^2 \dx,
\]
from which we can obtain
\[
    \frac{1}{2} \frac{\rm d}{\dt} \|\nabla \vv\|^2_2 + C \int_\Omega (1+ |\DD\vv|^2)^\frac{p(c)-2}{2}|\nabla \DD \vv|^2 \dx \le  C\|\fff\|_2^2 + C \|\eta\|_2^2 + C\|\nabla \vv\|_2^2\|\eta\|_2^2.
\]
Note further that \eqref{best1} implies $\eta \in L^2(Q_T)$. Thus by Gr\"onwall's inequality, we finally get the following uniform estimate
\begin{equation} \label{DDvest2}
    \sup_{0 \le t \le T} \| \nabla \vv\|_2^2 + \int_0^T \int_\Omega (1+ |\DD\vv|^2)^\frac{p(c)-2}{2}|\nabla \DD \vv|^2 \dx \dt \le C.
\end{equation}
Since $|\nabla^2\vv|\leq 3|\nabla\DD\vv|$, we can conclude that $\vv \in L^\infty(0,T;W^{1,2}(\Omega)) \cap  L^2(0,T;W^{2,2}(\Omega))$.  Furthermore, we shall insert \eqref{2DetaGNI} to \eqref{2Dgradeta} and apply Young's inequality to get
\[
    \|\nabla \eta\|_2^2 \le C \|(\overline{\DD}\vv)^\frac{p(c)-2}{2} \nabla \DD\vv\|_2^2 + C \|\eta\|_2^2\|\DD\vv\|_2^2 + C \|\eta\|_2^2.
\]
This, together with \eqref{DDvest2} and \eqref{best1}, implies
\begin{equation} \label{2DgradetaL2}
    \nabla \eta \in L^2(Q_T),
\end{equation}
which will be used in the later analysis.

\subsection{Higher integrability of $c_t$} \label{hct2d}
When we differentiated the stress tensor by the spatial variable, there appeared the $\nabla c$ term. Similarly, if we consider the time differentiation of the stress tensor, then the $\pa_t c$ term necessarily arises. Thus we need to control this term and we shall address it in Section \ref{hct2d}. In the two-dimensional domain, we know that
\[
    \vv \in L^\infty(0,T;W^{1,2}(\Omega)) \hookrightarrow  L^\infty(0,T;L^{s}(\Omega)) \quad \text{for} \quad 1 \le s < \infty.
\]
Furthermore, since we showed $\nabla c \in L^\infty(0,T;L^q(\Omega)) $ in \eqref{high_int_c}, we obtain
\[
\vv \cdot \nabla c \in L^\delta(Q_T)
\]
for some $\delta\in(4,q)$.
Moreover, from the assumption $ \ffg \in L^q(0,T;W^{1,q}(\Omega))$ and $c_0 \in W^{2,q}(\Omega)$, we can also see that ${\rm div} \, \ffg \in L^\delta(Q_T)$ and $\Delta c_0 \in L^\delta(\Omega)$.
Now for the concentration equation, we may write
\[
    \pa_t c - \Delta c =  h,
\]
where $h = {\rm div} \, (-c \vv - \ffg)$. Then by Lemma \ref{MR}, we have
\[
    \|\pa_t c\|_{L^\delta(Q_T)} + \|\Delta c\|_{L^\delta(Q_T)} \le C \| \vv \cdot \nabla c + {\rm div} \, \ffg \|_{L^\delta(Q_T)} + C \|\Delta c_0\|_{L^\delta(\Omega)}.
\]
Therefore, we can conclude for some $\delta\in(4,q)$ that
\begin{equation} \label{ct}
    \pa_t c \in L^\delta(Q_T)\quad{\rm{and}}\quad\Delta c \in L^\delta(Q_T).
\end{equation}

\subsection{$L^2$-estimate of $\partial_t \vv$} \label{2Dct}
Let us now investigate the regularity of $\pa_t \vv$. If we multiply $\pa_t \alpha_i^n$ to \eqref{aux11} and sum over $i=1,\ldots,n$, then we have
\begin{equation} \label{timetest}
    \int_\Omega \partial_t \vv \cdot \partial_t \vv \dx - \int_\Omega (\vv \otimes \vv) : \nabla (\partial_t \vv) \dx + \int_\Omega S(c,\DD \vv) : \DD (\partial_t \vv) \dx = \int_\Omega \fff \cdot \partial_t \vv \dx.
\end{equation}
For the stress tensor term, note that
\[
\frac{\rm d}{\dt} \bigg( (1+|\DD \vv|^2 )^\frac{p(c)}{2} \frac{1}{p(c)} \bigg) = \frac{\rm d}{\dt} \bigg((1+ |\DD \vv|^2)^\frac{p(c)}{2}\bigg) \frac{1}{p(c)} - \frac{ p'(c)}{p^2(c)}( 1+ |\DD \vv|^2)^\frac{p(c)}{2}  \partial_t c.
\]
Furthermore, there holds
\[
\frac{\rm d}{\dt} \bigg((1+ |\DD \vv|^2)^\frac{p(c)}{2}\bigg) \frac{1}{p(c)} = (1+ |\DD \vv|^2)^\frac{p(c)-2}{2} \DD \vv : \DD (\partial_t \vv) + \frac{1}{2} \frac{ p'(c)}{p(c)}(1+ |\DD\vv|^2)^\frac{p(c)}{2}  \partial_t c \log(1+ |\DD \vv|^2),
\]
which leads us to
\begin{equation} \label{vtest3}
\begin{aligned}
    (1+|\DD \vv|^2)^\frac{p(c)-2}{2} \DD \vv : \DD (\partial_t \vv)
    &=\frac{\rm d}{\dt}  \bigg( (1+|\DD \vv|^2 )^\frac{p(c)}{2} \frac{1}{p(c)} \bigg) + \frac{p'(c)}{p^2(c)}( 1+ |\DD \vv|^2)^\frac{p(c)}{2}  \partial_t c\\
    &\hspace{5mm} -\frac{1}{2}\frac{p'(c)}{p(c)} (1+ |\DD\vv|^2)^\frac{p(c)}{2}  \partial_t c \log(1+ |\DD \vv|^2).
\end{aligned}
\end{equation}
Therefore, from \eqref{timetest} with the assumption $p \in W^{1,\infty}(\mathbb{R})$, we obtain
\begin{equation}\label{3_5_5}
\begin{aligned}
    C\|\partial_t \vv\|_2^2 +  \frac{\rm d}{\dt} \int_\Omega (1+|\DD \vv|^2)^\frac{p(c)}{2} \frac{1}{p(c)} \dx &\leq C\int_\Omega (1+|\DD\vv|^2)^\frac{p(c)}{2} \partial_t c \dx+ C\int_\Omega (\vv \otimes \vv): \nabla \partial_t \vv \dx\\
    &\hspace{5mm}+ C\int_\Omega(1+|\DD\vv|^2)^\frac{p(c)}{2} \log(2+|\DD \vv|) \partial_t c \dx \\
    &\hspace{5mm}+ C\int_\Omega \fff \cdot \partial_t \vv \dx \\
    & =: {\rm{L}}_1 + {\rm{L}}_2 + {\rm{L}}_3 + {\rm{L}}_4.
\end{aligned}
\end{equation}
 We start with the estimate for ${\rm{L}}_1$. Recalling the definition of $\eta$ from \eqref{eta}, we have by Young's inequality that
\[
{\rm{L}}_1 \le C \| \eta\|_{4}^{4} + C \|\partial_t c \|_{2}^{2}.
\]
Now we recall \eqref{etabetaest} which implies that
\begin{equation} \label{eta44est}
\begin{split} 
\|\eta\|_4^4 & \le C \bigg( \|(\overline{\DD} \vv)^\frac{p(c)-2}{2} \nabla \DD \vv\|_2^2\|\eta\|_2^2 + \|\log(2+ |\DD \vv|) \|_\alpha^4 \|\eta\|_2^4 \bigg)  \\
& \le  C \bigg( \|(\overline{\DD} \vv)^\frac{p(c)-2}{2} \nabla \DD \vv\|_2^2 + \| 1+ |\DD \vv| \|_2^\frac{8}{\alpha} \|\eta\|_2^2 \bigg) \|\eta\|_2^2,
\end{split}
\end{equation}
where $\alpha$ is given in \eqref{alphaq2d}. Consequently, we have the following estimate for ${\rm{L}}_{1}$:
\begin{equation}\label{3_5_1}
	{\rm{L}}_1 \le C \bigg( \|(\overline{\DD} \vv)^\frac{p(c)-2}{2} \nabla \DD \vv\|_2^2 + \| 1+ |\DD \vv| \|_2^\frac{8}{\alpha} \|\eta\|_2^2 \bigg) \|\eta\|_2^2 + C \|\partial_t c \|_2^2.
\end{equation}
Secondly, for the convective term, a simple application of H\"older's inequality together with Young's inequality gives 
\begin{equation}\label{3_5_2}
    {\rm{L}}_2 = -C \int_\Omega \nabla \vv : ( \vv \otimes \partial_t\vv) \dx \le C \int_\Omega |\nabla \vv | |\vv| |\partial_t\vv| \dx \le C \int_\Omega |\nabla \vv |^2|\vv|^2 \dx + \epsilon \| \partial_t \vv\|_2^2.
\end{equation}
For the third term, we estimate ${\rm{L}}_3$ by H\"older's inequality and Young's inequality as
\begin{equation}\label{3_5_3}
{\rm{L}}_3 \le C \|\eta\|_4^4 + C \|\log(2+ |\DD\vv|)\|_\gamma^\gamma + C\|\partial_t c\|_\delta^\delta \le C\|\eta\|_4^4 + C \| 1+ |\DD \vv|\|_2^2 + C \|\partial_t c \|_\delta^\delta ,
\end{equation}
where $1 = \frac{1}{2} + \frac{1}{\gamma} + \frac{1}{\delta}$ with $\delta$ from \eqref{ct} and $\gamma$ large enough.
Finally, we can see from Young's inequality that
\begin{equation}\label{3_5_4}
{\rm{L}}_4 \le C \| \fff\|_2^2 + \epsilon \|\partial_t \vv\|_2^2.
\end{equation}
By collecting all the estimates \eqref{3_5_1}-\eqref{3_5_4} and inserting them into \eqref{3_5_5}, we have
\begin{equation}\label{vtest45}
\begin{aligned}
    & C\|\partial_t \vv\|_{L^2(\Omega)}^2 +  \frac{\rm d}{\dt} \int_\Omega (1+|\DD \vv|^2)^\frac{p(c)}{2} \frac{1}{p(c)} \dx \\
    & \hspace{4mm}\le C \bigg(\|(\overline{\DD} \vv)^\frac{p(c)-2}{2} \nabla \DD \vv\|_2^2 + \| 1+ |\DD \vv| \|_2^\frac{8}{\alpha} \|\eta\|_2^2 \bigg) \|\eta\|_2^2 \\
    & \hspace{8mm}+ C \bigg( \|\partial_t c\|_2^2 +  \|\partial_t c\|_\delta^\delta + \| 1+ |\DD \vv|\|_2^2 + \int_\Omega |\nabla \vv|^2 |\vv|^2 \dx + \|\fff\|_2^2 \bigg).
\end{aligned}
\end{equation}
First, note that 
\begin{equation} \label{eta22p}
\|\eta\|_2^2 = \int_\Omega (1+ |\DD \vv|^2)^\frac{p(c)}{2} \dx = \int_\Omega (1+ |\DD \vv|^2)^\frac{p(c)}{2} \frac{p(c)}{p(c)} \dx \le \int_\Omega (1+|\DD \vv|^2)^\frac{p(c)}{2} \frac{p^+}{p(c)} \dx.
\end{equation}
Also from \eqref{best1} and \eqref{DDvest2}, we have $\|\eta\|_2, \|(\overline{\DD} \vv)^\frac{p(c)-2}{2} \nabla \DD \vv\|_2 \in L^2(0,T)$ and $\|\DD \vv\|_2 \in L^\infty(0,T)$. Thus we can estimate the first term on the right-hand side of \eqref{vtest45} as 
\[
C\bigg(\|(\overline{\DD} \vv)^\frac{p(c)-2}{2} \nabla \DD \vv\|_2^2 + \| 1+ |\DD \vv| \|_2^\frac{8}{\alpha} \|\eta\|_2^2 \bigg) \|\eta\|_2^2\leq F(t) \int_\Omega (1+ |\DD \vv|^2)^\frac{p(c)}{2} \frac{1}{p(c)} \dx
\]
for some $F(t) \in L^1(0,T)$. Next, from \eqref{ct} and the assumption of the theorem, we see that 
\begin{equation*}
    \|\partial_t c\|_2^2, \,\,\|\partial_t c\|_\delta^\delta, \,\,\|\fff\|^2_2\in L^1(0,T).
\end{equation*} 
Furthermore, since
$L^2(0,T;W^{2,2}(\Omega)) \hookrightarrow L^2(0,T;L^\infty(\Omega))$, we have from \eqref{DDvest2} that
\begin{equation*}
    \int_0^T \int_\Omega |\nabla \vv|^2 |\vv|^2 \dx \dt \le \int^T_0 \left(\|\vv\|^2_{L^{\infty}(\Omega)}\int_{\Omega}\|\nabla \vv\|^2\dx \right)\dt\le C.
\end{equation*}
Thus, the second term on the right-hand side of \eqref{vtest45} can be bounded above by some $G(t)\in L^1(0,T)$. Therefore we can write \eqref{vtest45} as
\[
C\|\partial_t \vv\|_{L^2(\Omega)}^2 +  \frac{\rm d}{\dt} \int_\Omega (1+|\DD \vv|^2)^\frac{p(c)}{2} \frac{1}{p(c)} \dx \le F(t) \int_\Omega (1+|\DD \vv|^2)^\frac{p(c)}{2} \frac{1}{p(c)} \dx + G(t)
\]
and if we use Gr\"onwall's inequality with the assumption $\vv_0 \in W^{1,p^+}(\Omega)$, we can obtain the desired regularity
\[
\partial_t \vv \in L^2(0,T;L^2(\Omega)).
\]
Furthermore, we also have shown that
\begin{equation} \label{etainfty2D}
\int_\Omega (1+ |\DD \vv|^2)^\frac{p(c)}{2} \dx \in L^\infty(0,T).
\end{equation}
From the above together with \eqref{eta44est}, we conclude
\begin{equation} \label{eta44L1}
    \|\eta\|_4^4 \in L^1(0,T),
\end{equation}
which will be used in the later analysis.

\section{Uniform estimates in $3$D}

In this section, we will derive uniform estimates in the three-dimensional domain. Note that the energy estimates \eqref{best1} and \eqref{best2} in Section \ref{best} still hold for $d=3$. Also, provided that $p^-\geq\frac{d+2}{2}$, the estimate \eqref{high_int_c} concerning the integrability of $\nabla c$ in Section \ref{hgradc} holds as well. In the following subsections, we aim to obtain uniform estimates in the three-dimensional domain whose derivations are different from the ones for the two-dimensional case.

\subsection{$W^{2,2}$-estimate of $\vv$}

As we did in the 2D estimates, we may write
\[
\int_{\Omega} \pa_t \vv \cdot (-\Delta \vv) \dx - \int_{\Omega} (\vv\otimes\vv) : \nabla (-\Delta \vv) \dx + \int_{\Omega} \SSS(c,\DD \vv):\DD (-\Delta \vv) \dx = \int_{\Omega} \fff \cdot (-\Delta \vv) \dx,
\]
which can be rewritten as
\begin{equation} \label{W22est3d2}
\frac{1}{2}\frac{\rm d}{{\rm d}t}\int_{\Omega} |\nabla \vv|^2 \dx  + {\rm I_2} = {\rm I_1} + \int_{\Omega} \fff \cdot (-\Delta \vv) \dx,
\end{equation}
with
\[
{\rm I_1}=\int_{\Omega} (\vv\otimes\vv) : \nabla (-\Delta \vv) \dx\quad{\rm and}\quad{\rm I_2}=\int_{\Omega} \SSS(c,\DD \vv):\DD (-\Delta \vv) \dx.
\]
From the repetitive uses of the integration by parts, we have
\begin{equation} \label{gradv3}
 {\rm I_1} \le \int_\Omega |\nabla \vv|^3 \dx
\end{equation}
and the application of the chain rule yields
\[
   {\rm I_2} = \int_{\Omega}  \frac{\pa \SSS(c,\DD)}{\pa \DD}:(\pa_{x_k} \DD\vv\otimes\pa_{x_k}\DD\vv)\dx +\int_{\Omega} \pa_{x_k} c \left( \frac{\pa\SSS(c,\DD)}{\pa c} : \pa_{x_k}\DD\vv \right)
  \dx 
=: {\rm E_1} + {\rm E_2}.
\]
As before, it is straightforward to see that
\begin{equation} \label{E13d}
     {\rm E_1} \ge C \int_\Omega ( 1+ |\DD \vv|^2)^\frac{p(c)-2}{2} |\nabla \DD \vv|^2 \dx.
\end{equation}
For the second term, we have as before that
\begin{equation} \label{E23d}
    |{\rm E_2}| \le C \|\nabla c\|_q \| \log(2+ |\DD \vv|) \|_\alpha \| (1+ |\DD \vv|^2)^\frac{p(c)}{4} \|_3 \|  ( 1+ |\DD \vv|^2)^\frac{p(c)-2}{4} \nabla \DD \vv\|_2,
\end{equation}
where
\begin{equation} \label{qalpha3d}
    1= \frac{1}{q} + \frac{1}{\alpha}+\frac{1}{3} + \frac{1}{2}
\end{equation}
with $q>6$ as given in Theorem \ref{maintheorem} and large enough $\alpha>6$. We shall estimate ${\rm E_2}$ in a similar way as we did in the 2D case. First of all, from \eqref{high_int_c}, we have $\sup_{0\le t \le T}\|\nabla c(t)\|_q \le C$. Next, by recalling the definition \eqref{eta}, we can rewrite \eqref{E23d} as
\begin{equation} \label{E23d2}
|{\rm E_2}| \le C \|\log(2+ |\DD\vv|)\|_\alpha \|\eta\|_3 \|(\overline{\DD}\vv)^\frac{p(c)-2}{2} \nabla \DD\vv\|_2.
\end{equation}
By the interpolation inequality, there holds
\begin{equation} \label{3DetaGNI}
\|\eta\|_3 \le C \|\eta\|_2^\frac{1}{2} \|\nabla \eta\|_2^\frac{1}{2}.
\end{equation}
Note that by the assumption $p \in W^{1,\infty}(\mathbb{R})$, H\"older's inequality and \eqref{high_int_c}, we get
\begin{equation} \label{3Dgradeta}
\|\nabla \eta\|_2 \le C \left( \|(\overline{\DD}\vv)^\frac{p(c)-2}{2} \nabla \DD\vv\|_2 + \|\eta\|_3 \|\log(2+ |\DD\vv|)\|_\alpha\right)\hspace{-0.1cm}.
\end{equation}
Therefore,  together with Young's inequality, we obtain
\begin{equation} \label{etabeta3d}
    \|\eta\|_3 \le C \left(\|(\overline{\DD}\vv)^\frac{p(c)-2}{2} \nabla \DD\vv\|_2^\frac{1}{2}\|\eta\|_2^\frac{1}{2} + \|\log(2+ |\DD\vv|)\|_\alpha \|\eta\|_2\right)\hspace{-0.1cm}.
\end{equation}
 Now if we insert \eqref{etabeta3d} to \eqref{E23d2} and use  Young's inequality, we have
\[ 
        |{\rm E_2}| \le \epsilon \|(\overline{\DD}\vv)^\frac{p(c)-2}{2} \nabla \DD \vv\|_2^2 + C \|\log(2+ |\DD \vv|)\|_\alpha^4 \|\eta\|^2_2.
\]
Since $\log(2+ |\DD\vv|) \le C (1+ |\DD\vv|)^\frac{2}{\alpha}$, we get $\|\log(2+ |\DD \vv|)\|_\alpha^4 \le C + C\|\DD\vv\|_2^\frac{8}{\alpha}$. As we know $\frac{8}{\alpha} <2 $, we can apply Young's inequality to get 
\[
    \|\log(2+|\DD\vv|)\|_\alpha^4 \le C + C \|\DD \vv\|_2^2.
\]
Thus we have
\begin{equation} \label{E23d_2}
    |{\rm E_2}| \le \epsilon \|(\overline{\DD}\vv)^\frac{p(c)-2}{2} \nabla \DD \vv\|_2^2 + C \|\eta\|_2^2 + C \|\DD\vv\|_2^2 \|\eta\|_2^2.
\end{equation}
Now if we insert \eqref{gradv3}, \eqref{E13d} and \eqref{E23d_2} to \eqref{W22est3d2}, and apply Young's inequality, then we have
\begin{equation} \label{W22est3d3}
       \frac{1}{2} \frac{\rm d}{\dt} \|\nabla \vv\|^2_2 + C\int_\Omega (1+ |\DD \vv|^2)^{\frac{p(c)-2}{2}} |\nabla \DD \vv|^2 \dx 
       \le  C\|\fff\|_2^2 + C \|\eta\|_2^2 + C\|\nabla \vv\|_2^2\|\eta\|_2^2 + \|\nabla \vv\|_3^3 .
\end{equation}
Now, the only difference compared to the 2D case is the presence of $\|\nabla \vv\|_3^3$ term. We will handle this term by the similar argument used in \cite{malek}. For the case of $p^- \ge 3$, as we have \eqref{best1}, we obtain $\nabla \vv \in L^3(Q_T)$. Then we can apply Gr\"onwall's inequality and obtain the desired result immediately.
Thus we may assume $p^-<3$ and let us write
\begin{equation} \label{grad33est1}
    \|\nabla \vv\|_3^3 = \|\nabla \vv\|_3^{3(1-\tau) + 3 \tau}
\end{equation}
for some $\tau\in(0,1)$. Note by Korn's inequality and the Sobolev inequality that
\[
    \int_\Omega |\nabla \vv|^{3p^-} \dx \le C \int_\Omega \big(|\DD \vv|^{\frac{p^-}{2}} \big)^6 \dx \le C \bigg( \int_\Omega \big| \nabla \big( |\DD\vv|^\frac{p^-}{2}\big) \big|^2 \dx \bigg)^3,
\]
which implies that
\begin{equation} \label{grad3pest1}
    \|\nabla \vv\|_{3p^-}^{p^-} \le C \int_\Omega (1+ |\DD \vv|^2)^\frac{p^--2}{2} |\nabla \DD \vv|^2 \dx \le C \int_\Omega (1+ |\DD\vv|^2)^\frac{p(c)-2}{2} |\nabla \DD\vv|^2 \dx.
\end{equation}
Now by the interpolation inequality, we have
\begin{align*}
    & \|\nabla \vv\|_3^{3(1-\tau)} \le \|\nabla \vv\|_2^{3(1-\tau)\frac{2p^--2}{3p^--2}} \|\nabla \vv\|_{3p^-}^{(1-\tau)\frac{3p^-}{3p^--2}}, \\
    & \|\nabla \vv\|_3^{3\tau} \le \|\nabla \vv \|_{p^-}^{3\tau \frac{p^- -1}{2}} \|\nabla \vv\|_{3p^-}^{3\tau\frac{3-p^-}{2}}.
\end{align*}
Hence, from \eqref{grad33est1}, we can write 
\[
    \|\nabla \vv\|_3^3 \le \|\nabla \vv\|_2^{r_1} \|\nabla \vv\|_{p^-}^{r_2} \|\nabla \vv\|_{3p^-}^{r_3},
\]
where
\[
    r_1 = 3(1-\tau)\frac{2p^--2}{3p^--2},\quad
    r_2 = 3\tau\frac{p^--1}{2}\quad{\rm{and}}\quad
    r_3 = (1-\tau)\frac{3p^-}{3p^--2} + 3\tau \frac{3-p^-}{2}.
\]
Now applying Young's inequality with \eqref{grad3pest1} yields
\begin{equation} \label{grad33est2}
    \|\nabla \vv\|_3^3 \le C \left[(\| \nabla \vv\|_2^2)^\frac{r_1}{2} ( 1+ \|\nabla \vv\|_{p^-})^{r_2} \right]^{s'} + \epsilon \int_\Omega (1+ |\DD \vv|^2)^\frac{p(c)-2}{2} |\nabla \DD \vv|^2 \dx,
\end{equation}
where $s = \frac{p^-}{r_3}$. We aim to choose $\tau$ satisfying $r_2 s' = p^-$, which motivates us to write
\[
    1= \frac{1}{s} + \frac{1}{s'} = \frac{r_3}{p^-} + \frac{r_2}{p^-} = \frac{3\tau}{p^-} + \frac{3(1-\tau)}{3p^--2}.
\]
Then it follows that
\[
    \tau = \frac{p^-(3p^--5)}{6p^--6}, \quad 1-\tau = \frac{(3-p^-)(3p^--2)}{6p^--6} \quad \text{and} \quad s' = \frac{4}{3p^--5}.
\]
Since $\frac{5}{2} \le p^- <3$, we can confirm that $0<\tau<1$. Hence by putting \eqref{grad33est2} to \eqref{W22est3d3}, we obtain 
\begin{equation*}
     \frac{1}{2} \frac{\rm d}{\dt} \|\nabla \vv\|^2_2 + C \int_\Omega  (1+|\DD\vv|^2)^\frac{p(c)-2}{2}|\nabla \DD \vv|^2 \dx
     \le  C\|\fff\|_2^2 + C \|\eta\|_2^2 + C\|\eta\|_2^2 \|\nabla \vv\|_2^2 + C (1+ \|\nabla \vv\|_{p^-}^{p^-}) (\|\nabla \vv\|_2^2)^\frac{r_1s'}{2}.
\end{equation*}
Now note that $\frac{r_1s'}{2} \le 1$ is equivalent to $\frac{11}{5} \le p^-$ which is satisfied from our assumption. Thus we apply Gr\"onwall's inequality to get 
\begin{equation} \label{W22est3dlast}
    \sup_{0 \le t \le T} \| \nabla \vv\|_2^2 + \int_0^{T}  \int_\Omega (1+|\DD\vv|^2)^\frac{p(c)-2}{2} |\nabla \DD \vv|^2 \dx \dt \le C.
\end{equation}
Therefore, we have obtained the desired regularity $\vv \in L^\infty(0,T;W^{1,2}(\Omega)) \cap L^2(0,T;W^{2,2}(\Omega))$ in the three-dimensional domain. Furthermore, by inserting \eqref{3DetaGNI} to \eqref{3Dgradeta} and applying Young's inequality, we have 
\[
    \|\nabla \eta\|_2^2 \le C \|(\overline{\DD}\vv)^\frac{p(c)-2}{2} \nabla \DD\vv\|_2^2 + C \|\eta\|_2^2\|\DD\vv\|_2^2 + C \|\eta\|_2^2,
\]
which together with \eqref{W22est3dlast} and \eqref{best1} implies
\begin{equation} \label{3DgradetaL2}
    \nabla \eta \in L^2(Q_T).
\end{equation}

\subsection{Higher integrability of $c_t$} \label{hct3d}
We first note that
\[
 \vv \in L^\infty(0,T;W^{1,2}(\Omega)) \hookrightarrow L^\infty(0,T;L^6(\Omega))
\]
and $\nabla c \in L^\infty(0,T;L^q(\Omega))$ from \eqref{high_int_c}, which lead us to have
\begin{equation*}
    \vv \cdot \nabla c \in L^\delta(Q_T)
\end{equation*}
with $\delta\in(3,3+\epsilon)$ for some $\epsilon>0$.
As we did before, from Lemma \ref{MR}, we get 
\begin{equation} \label{ct3d}
    \pa_t c \in L^\delta(Q_T) \quad \text{and} \quad
    \Delta c \in L^\delta(Q_T).
\end{equation}

\subsection{$L^2$ estimate of $\pa_t \vv$}
As we did before in \eqref{timetest}-\eqref{3_5_5}, we have 
\begin{equation} \label{L2patv3D}
\begin{aligned}
    C\|\partial_t \vv\|_2^2 +   \frac{\rm d}{\dt} \int_\Omega (1+|\DD \vv|^2)^\frac{p(c)}{2} \frac{1}{p(c)} \dx 
    &\le C\int_\Omega (1+|\DD\vv|^2)^\frac{p(c)}{2} \partial_t c \dx+ C\int_\Omega (\vv \otimes \vv): \nabla \partial_t \vv \dx\\    &\hspace{4mm}+C\int_\Omega(1+|\DD\vv|^2)^\frac{p(c)}{2} \log(1+|\DD \vv|)\partial_t c  \dx\\
    &\hspace{4mm} + C\int_\Omega \fff \cdot \partial_t \vv \dx\\
    & =: {\rm L_1} + {\rm L_2} + {\rm L_3} + {\rm L_4}. 
\end{aligned}
\end{equation}
First, by Young's inequality, we have
\begin{equation*}
    {\rm L_1} = C \int_\Omega \eta^2 \pa_t c \dx \le C\|\eta\|_{3}^{3} + C \|\pa_t c \|_{3}^{3},
\end{equation*}
where $\eta$ is defined in \eqref{eta}. Now we recall \eqref{etabeta3d}, which implies that
\begin{equation*}
    \|\eta\|_3^3 \le C \|(\overline{\DD}\vv)^\frac{p(c)-2}{2} \nabla \DD \vv\|_2^\frac{3}{2} \|\eta\|_2^\frac{3}{2} + C\|\log(2+ |\DD \vv|)\|_\alpha^3 \|\eta\|_2^3,
\end{equation*}
where $\alpha$ is given in \eqref{qalpha3d}. Then the application of Young's inequality yields
\begin{equation} \label{eta33est}
\begin{split}
    \|\eta\|_3^3 &\le C \bigg(\|(\overline{\DD} \vv)^\frac{p(c)-2}{2} \nabla \DD \vv\|_2^2 \|\eta\|_2^2 + \|\log(2+ |\DD\vv|)\|_\alpha^4\|\eta\|_2^4  \bigg) +C \\
    &\le C \bigg( \|(\overline{\DD}\vv)^\frac{p(c)-2}{2} \nabla \DD \vv\|_2^2 + \|1+|\DD\vv|\|_2^\frac{8}{\alpha} \|\eta\|_2^2 \bigg) \|\eta\|_2^2 + C.
\end{split}
\end{equation}
Thus we get
\begin{equation}  \label{L_13D}
    {\rm L_1} \le C\|\pa_t c \|_3^3 + C \bigg( \|(\overline{\DD}\vv)^\frac{p(c)-2}{2} \nabla \DD \vv\|_2^2  +  \|1+|\DD\vv|\|_2^\frac{8}{\alpha} \|\eta\|_2^2 \bigg) \|\eta\|_2^2 +C .
\end{equation}
For the $L_2$ term, since $\vv \in L^\infty(0,T;L^6(\Omega))$, we obtain  by H\"older's inequality and Young's inequality that
\begin{equation} \label{L_23D}
\begin{aligned}
    {\rm L_2} 
    &= -C\int_\Omega \nabla \vv : ( \vv \otimes \pa_t \vv) \dx \le C \int_\Omega |\nabla \vv|^2 |\vv|^2 \dx + \epsilon \|\pa_t \vv\|_2^2 \\
    &\le C\|\nabla \vv\|_3^2 \|\vv\|_6^2 + \epsilon\|\pa_t \vv\|_2^2 \le C\|\nabla \vv\|_3^2 + \epsilon\|\pa_t \vv\|_2^2.
\end{aligned}
\end{equation}
Next, by Young's inequality, we obtain that
\begin{equation}  \label{L_33D}
    {\rm L_3} \le C \|\eta\|_3^3 + C\|\log(2+ |\DD \vv|)\|_\gamma^\gamma + C\|\pa_t c\|_\delta^\delta \le C\|\eta\|_3^3 + C\|1+|\DD\vv|\|_2^2 + C\|\pa_t c\|_\delta^\delta,
\end{equation}
where $1= \frac{2}{3} + \frac{1}{\gamma} + \frac{1}{\delta}$ with arbitrarily large $\gamma$ and $\delta$ from \eqref{ct3d}. 
Finally, ${\rm L_4}$ can be easily estimated by Young's inequality as
\begin{equation}  \label{L_43D}
    {\rm L_4} \le C \|\fff\|_2^2 + \epsilon\|\pa_t \vv\|_2^2.
\end{equation}
Now we combine \eqref{L_13D}-\eqref{L_43D} and insert them into \eqref{L2patv3D} to get
\begin{equation} \label{vtest3dl}
\begin{aligned}
    C\|\partial_t &\vv\|_{L^2(\Omega)}^2 +  \frac{\rm d}{\dt} \int_\Omega (1+|\DD \vv|^2)^\frac{p(c)}{2} \frac{1}{p(c)} \dx  \\
    & \le C \bigg(\|(\overline{\DD} \vv)^\frac{p(c)-2}{2} \nabla \DD \vv\|_2^2 + \| 1+ |\DD \vv| \|_2^\frac{8}{\alpha} \|\eta\|_2^2 \bigg) \|\eta\|_2^2  \\
    & \hspace{4mm}+ C \bigg( \|\partial_t c\|_3^3 +  \|\partial_t c\|_\delta^\delta + \| 1+ |\DD \vv|\|_2^2 + \|\nabla \vv\|_3^2 + \|\fff\|_2^2 +1 \bigg)\hspace{-0.05cm}. 
\end{aligned}
\end{equation}
As we did in the 2D case, \eqref{best1} and \eqref{W22est3dlast} imply $\|\eta\|_2, \|(\overline{\DD}\vv)^\frac{p(c)-2}{2} \nabla \DD \vv\|_2 \in L^2(0,T)$ and $\|\DD \vv\|_2 \in L^\infty(0,T)$. Together with \eqref{eta22p}, the first term on the right-hand side of \eqref{vtest3dl} can be written as
\[
C\bigg(\|(\overline{\DD} \vv)^\frac{p(c)-2}{2} \nabla \DD \vv\|_2^2 + \| 1+ |\DD \vv| \|_2^\frac{8}{\alpha} \|\eta\|_2^2 \bigg) \|\eta\|_2^2 \le  F(t) \int_\Omega(1+|\DD\vv|^2)^\frac{p(c)}{2} \frac{1}{p(c)} \dx
\]
for some $F(t) \in L^1(0,T)$. Also by using the facts $\vv \in L^2(0,T;W^{2,2}(\Omega)) \hookrightarrow L^2(0,T;W^{1,6}(\Omega))$, $\fff \in L^2(Q_T)$ and \eqref{ct3d}, we finally have
\begin{equation*}
    C\|\partial_t \vv\|_{L^2(\Omega)}^2 +  \frac{\rm d}{\dt} \int_\Omega (1+|\DD \vv|^2)^\frac{p(c)}{2} \frac{1}{p(c)} \dx \le F(t) \int_\Omega (1+|\DD \vv|^2)^\frac{p(c)}{2} \frac{1}{p(c)} \dx + G(t)
\end{equation*}
for some $ G(t) \in L^1(0,T)$. Hence by Gr\"onwall's inequality with the assumption $\vv_0 \in W^{1,p^+}(\Omega)$, we conclude that
\begin{equation*}
    \pa_t \vv \in L^2(0,T;L^2(\Omega)).
\end{equation*}
Furthermore, we can also obtain
\begin{equation} \label{3DDVLinfty}
    \int_\Omega (1 + |\DD\vv|^2)^\frac{p(c)}{2}  \dx \in L^\infty(0,T),
\end{equation}
which together with \eqref{eta33est} implies that
\begin{equation} \label{eta33L1}
    \|\eta\|_3^3 \in L^1(0,T).
\end{equation}

\section{Proof of Theorem \ref{maintheorem}}
From the uniform estimates derived in the previous sections, in both 2D and 3D cases, we can extract a (not relabeled) subsequence with respect to $n\in\mathbb{N}$ such that
\[
    \begin{aligned}
        \vv^n &\rightharpoonup \vv &&\mbox{weakly in} \quad L^2(0,T;W^{2,2}(\Omega)), \\
     \vv^n & \overset{\ast} {\rightharpoonup} \vv && \mbox{weakly-}^*\,\,{\rm{in}} \quad L^\infty(0,T;W^{1,2}(\Omega)), \\
     \partial_t \vv^n &\rightharpoonup \partial_t \vv && \mbox{weakly in} \quad L^2(0,T;L^2(\Omega)), \\
     c^n & \rightharpoonup c && \text{weakly in} \quad L^2(0,T;W^{1,2}(\Omega)), \\
     c^n & \overset{\ast}  {\rightharpoonup} c && \text{weakly-}^*\,\,{\rm{in}} \quad L^\infty(0,T;W^{1,q}(\Omega)), \\
     \pa_t c^n & \rightharpoonup \pa_t c && \text{weakly in} \quad L^2(0,T;L^2(\Omega)). 
    \end{aligned}
\]
By applying the Aubin--Lions lemma presented in Lemma \ref{ALL}, we have that 
\begin{equation} \label{conv_v}
    \vv^n \to \vv \quad \mbox{strongly in} \quad L^2(0,T;W^{1,2}(\Omega)),
\end{equation}
from which we can obtain
\begin{equation*}
    \DD\vv^n \to \DD\vv \quad \text{a.e. in} \quad Q_T.
\end{equation*}
For the concentration, by the Aubin--Lions lemma again, we have
\begin{equation*}
    c^n \to c \quad \text{strongly in} \quad L^2(0,T;L^2(\Omega)),
\end{equation*}
which implies that
\begin{equation*}
    c^n \to c \quad \text{a.e. in} \quad Q_T.
\end{equation*}
Therefore, by the continuity of $\SSS(\cdot,\cdot)$, we can see that
\begin{equation*}
    \SSS(c^n,\DD\vv^n) \to \SSS(c,\DD\vv) \quad \text{a.e. in} \,\,\,Q_T.
\end{equation*}
Moreover, we have from \eqref{best1} that
\begin{equation*}
    \|\SSS(c^n,\DD\vv^n)\|_{L^1(Q_T)} \le C.
\end{equation*}
Thus by Vitali's theorem, we obtain the identification
\begin{equation*}
    \SSS(c^n,\DD\vv^n) \to \SSS(c,\DD\vv) \quad \text{a.e. in} \quad L^1(Q_T).
\end{equation*}
Also, as a consequence of \eqref{conv_v} we have
\begin{equation*}
    (\vv^n \cdot \nabla) \vv^n \to (\vv \cdot \nabla) \vv \quad \mbox{in} \quad L^\frac{4}{3}(Q_T).
\end{equation*}
Now we shall temporarily choose $a\in(1,2)$ which will be specified later. We first note that 
\begin{equation}\label{nablaS}
\begin{split}
    |\nabla \SSS(c^n,\DD\vv^n)|^a &\le C(1+|\DD\vv^n|^2)^{a\frac{p(c^n)-2}{2}} |\nabla \DD \vv^n|^a  \\ 
    & \hspace{4mm} + C (1+ |\DD\vv^n|^2)^{a\frac{p(c^n)-1}{2}}|\log(2+|\DD\vv^n|)|^a |\nabla c^n|^a  \\
    &=: {\rm M_1} + {\rm M_2}. 
\end{split}
\end{equation}
Let us consider the cases of 2D and 3D separately. First of all, in the case of 2D, we set $a=\frac{4}{3}$. Then by Young's inequality, ${\rm M_1}$ can be estimated as
\begin{align}
    {\rm M_1} &= C (1+ |\DD\vv^n|^2)^{\frac{4}{3}\frac{p(c^n)-2}{4}}(1+|\DD\vv^n|^2)^{\frac{4}{3}\frac{p(c^n)-2}{4}} |\nabla \DD\vv^n|^\frac{4}{3} \nonumber \\ 
    & \le C (1+ |\DD\vv^n|^2)^{4\frac{p(c^n)-2}{4}} + C ( 1+ |\DD\vv^n|^2)^\frac{p(c^n)-2}{2} |\nabla \DD\vv^n|^2. \label{M_1est2}
\end{align}
From \eqref{DDvest2}, the second term in \eqref{M_1est2} belongs to $L^1(Q_T)$. For the first term, by recalling \eqref{eta} and \eqref{eta44L1}, we have
\begin{equation} \label{M1etaest2d2}
    (1+|\DD\vv^n|^2)^{4\frac{p(c^n)-2}{4}} \le (1+ |\DD\vv^n|^2)^{4\frac{p(c^n)}{4}} = \eta^4 \in L^1(Q_T).
\end{equation}
Thus we get ${\rm M_1} \in L^1(Q_T)$. For the estimate of ${\rm M_2}$, by Young's inequality
\[
{\rm M_2} \le C(1+|\DD\vv^n|^2)^{4{\frac{p(c^n)-1}{4}}} + |\log(2+ |\DD\vv^n|)|^\frac{4\alpha}{3} + |\nabla c^n|^q,
\]
where $1=\frac{2}{3} + \frac{1}{\alpha} + \frac{4}{3q}$ with $q>4$ as given in theorem \ref{maintheorem} and arbitrarily large $\alpha>0$. Also, it is straightforward that
\[
|\log(2+ |\DD\vv^n|)|^\frac{4\alpha}{3} \le C + C|\DD\vv^n|^2.
\]
Thus by \eqref{best1}, \eqref{high_int_c} and \eqref{M1etaest2d2}, we have ${\rm M_2} \in L^1(Q_T)$.

For the case of 3D, we let $a=\frac{6}{5}$ and then ${\rm M_1}$ can be estimated by Young's inequality as
\begin{align}
    {\rm M_1} &= C (1+ |\DD\vv^n|^2)^{\frac{6}{5}\frac{p(c^n)-2}{4}}(1+|\DD\vv^n|^2)^{\frac{6}{5}\frac{p(c^n)-2}{4}} |\nabla \DD\vv^n|^\frac{6}{5} \nonumber \\ 
    & \le C (1+ |\DD\vv^n|^2)^{3\frac{p(c^n)-2}{4}} + C ( 1+ |\DD\vv^n|^2)^\frac{p(c^n)-2}{2} |\nabla \DD\vv^n|^2. \label{M_1est3D}
\end{align}
From \eqref{W22est3dlast}, the second term in \eqref{M_1est3D} belongs to $L^1(Q_T)$. By recalling \eqref{eta} and \eqref{eta33L1}, we have
\begin{equation} \label{M1etaest3}
    (1+|\DD\vv^n|^2)^{3\frac{p(c^n)-2}{4}} \le (1+ |\DD\vv^n|^2)^{3\frac{p(c^n)}{4}} = \eta^3 \in L^1(Q_T),
\end{equation}
which implies that the first term of \eqref{M_1est3D} belongs to $L^1(Q_T)$. Hence we obtain ${\rm M_1} \in L^1(Q_T).$
For the estimate of ${\rm M_2}$, by Young's inequality
\[
    {\rm M_2} \le C(1+|\DD\vv^n|^2)^{3\frac{p(c^n)-1}{4}} + |\log(2+ |\DD\vv^n|)|^\frac{6\alpha}{5} + |\nabla c^n|^q,
\]
where $1= \frac{4}{5} + \frac{1}{\alpha} + \frac{6}{5q}$ with $q>6$ as given in theorem \ref{maintheorem} and large enough $\alpha>0$. The logarithmic term can be controlled as we did in the 2D case, and then by \eqref{best1}, \eqref{high_int_c} and \eqref{M1etaest3} we also have  ${\rm M_2} \in L^1(Q_T).$
Therefore, from \eqref{nablaS} with the above estimates, we have for both 2D and 3D cases that
\begin{equation*}
    \nabla \SSS(c^n,\DD\vv^n) \in L^\frac{6}{5} (Q_T).
\end{equation*}
Altogether, we obtain that
\begin{equation*}
    \|\pa_t \vv\|_{L^2(Q_T)} + \|\nabla\SSS(c,\DD\vv)\|_{L^\frac{6}{5}(Q_T)} + \| (\vv \cdot\nabla)\vv\|_{L^\frac{4}{3}(Q_T)} \le C.
\end{equation*}
Now by choosing $\ww \in \{\ww_i\}_{i=1}^\infty$ from \eqref{vbasis} and $\xi \in C_0^\infty([0,T])$, we deduce from \eqref{aux1} and the above inequality that
\begin{equation*}
    \int_{Q_T} \big(\pa_t \vv \cdot \ww \xi + (\vv \cdot \nabla) \vv \cdot \ww\xi  - {\rm div} \, (\SSS(c,\DD\vv)) \cdot \ww \xi \big)\dx\dt= \int_{Q_T} \fff \cdot \ww \xi \dx \dt.
\end{equation*}
By the density of smooth functions, we get
\begin{equation*}
    \int_{Q_T} \big(\pa_t \vv \cdot \ffpsi+ (\vv \cdot \nabla) \vv \cdot \ffpsi  - {\rm div} \, (\SSS(c,\DD\vv)) \cdot \ffpsi \big) \dx\dt = \int_{Q_T} \fff \cdot \ffpsi \dx \dt
\end{equation*}
for all $\ffpsi \in C_0^\infty([0,T];\mathcal{V})$. Now by De Rahm's theorem, there exists $\nabla \pi \in L^\frac{6}{5}(Q_T)$ (see \cite{breit, Simon1993}) and we finally get
\begin{equation*}
    \pa_t \vv + {\rm div} \, (\vv \otimes \vv) - {\rm div} \, (\SSS(c,\DD\vv)) + \nabla \pi = \fff \quad \text{a.e. in} \quad Q_T.
\end{equation*}
On the other hand, from \eqref{ct} and \eqref{ct3d}, we have for both 2D and 3D cases that
\begin{equation*}
    \Delta c \in L^2(Q_T).
\end{equation*}
Since $\vv \in L^2(0,T;L^\infty(\Omega))$, $c \in L^\infty(0,T;W^{1,q}(\Omega))$ and $q>2d$, it is straightforward to verify in both two and three-dimensional domains
\begin{equation*}
    \vv \cdot \nabla c \in L^2(Q_T).
\end{equation*}
In all, we have
\begin{equation*}
    \|\pa_t c\|_{L^2(Q_T)} + \|\, {\rm div} \, (c \vv)\|_{L^2(Q_T)} + \|\Delta c\|_{L^2(Q_T)} \le C.
\end{equation*}
Now we deduce from \eqref{aux33} and again by the density of smooth functions, we obtain
\begin{equation*}
    \int_{Q_T} \big(\pa_t c \,\phi + {\rm div}\,  (c\vv) \, \phi - \Delta c \, \phi \big)\dx\dt = \int_{Q_T} {\rm div}\, \ffg \, \phi \dx\dt
\end{equation*}
for all $\phi \in C_0^\infty([0,T];C^\infty(\Omega))$, which implies that
\begin{equation*}
    \pa_t c + {\rm div}\,(c\vv) - \Delta c = {\rm div} \, \ffg \quad \text{a.e. in} \quad Q_T.
\end{equation*}
Finally, we can see that $\vv, c \in C([0,T],L^2(\Omega))$ since $\vv, c \in L^2(0,T;W^{1,2}(\Omega))$ and $\pa_t \vv, \pa_t c \in L^2(0,T;L^2(\Omega))$. Moreover, it is obvious that $\|P^n \vv_0 - \vv_0\|_2 \to 0$ which completes the proof.

\begin{remark}
    Indeed, we have shown that $\nabla \SSS(c^n,\DD\vv^n) \in L^\frac{4}{3}(Q_T)$ in the two-dimensional domain. Therefore, by De Rahm's theorem, $\nabla \pi$ belongs to $L^\frac{4}{3}(Q_T)$ in 2D.
\end{remark}

\section{Proof of Theorem \ref{uniqueness theorem}}
In this section, let us prove Theorem \ref{uniqueness theorem} concerning the uniqueness of the strong solution constructed in Theorem \ref{maintheorem}. Let $(\vv_1,c_1)$ and $(\vv_2,c_2)$ be two solutions of the problem \eqref{main_sys}-\eqref{concentration}. First of all, we take the difference of \eqref{concentration} of $c_1$ and $c_2$ to get
\begin{equation} \label{c1c2equation}
    \pa_t(c_1-c_2) - \Delta(c_1-c_2) + {\rm div}\, (c_1\vv_1 - c_2 \vv_2) = 0.
\end{equation}
We then test $-\Delta (c_1-c_2)$ to the above equation and use Young's inequality to get
\begin{equation}\label{uniconcentration}
\begin{aligned}
    \frac{1}{2}\frac{\rm d}{\dt} \int_\Omega |\nabla (c_1-c_2)|^2 \dx + \int_\Omega |\Delta (c_1-c_2)|^2 \dx 
    &\le \int_\Omega|\vv_1 \cdot \nabla c_1 - \vv_2 \cdot \nabla c_2|^2 \dx\\
    &= \int_\Omega |(\vv_1-\vv_2) \cdot \nabla c_1  + \vv_2 \cdot (\nabla c_1 - \nabla c_2)|^2\dx\\
    &\le C \|\nabla c_1\|_\infty^2 \|\vv_1 - \vv_2 \|_2^2 + C \|\vv_2\|_\infty^2 \|\nabla (c_1 -c_2)\|_2^2.
\end{aligned}
\end{equation}
Secondly, take the difference of \eqref{main_sys} of $\vv_1$ and $\vv_2$ and test $\vv_1 - \vv_2$ to get
\begin{equation} \label{uniPQ}
    \frac{1}{2} \frac{\rm d}{\dt} \|\vv_1 - \vv_2\|_2^2 + {\rm P} = {\rm Q},
\end{equation}
where
\[
\begin{split}
    {\rm P} &:= \int_\Omega \left( \SSS(c_1,\DD\vv_1)-\SSS(c_2,\DD\vv_2)\right) : (\DD\vv_1-\DD\vv_2) \dx, \\
    {\rm Q} &:= - \int_\Omega \left( (\vv_1\cdot \nabla )\vv_1 - (\vv_2 \cdot \nabla ) \vv_2 \right) \cdot (\vv_1 - \vv_2) \dx.
\end{split}
\]
Note by H\"older's inequality and Young's inequality that
\[
\begin{split}
    {\rm Q}
    &= -\int_\Omega\left((\vv_1-\vv_2)\cdot \nabla \right)\vv_1 \cdot (\vv_1 - \vv_2) \dx \le \int_\Omega |\vv_1-\vv_2| |\nabla \vv_1||\vv_1 - \vv_2|  \dx \\
    &\le \|\vv_1 - \vv_2 \|_6 \|\nabla \vv_1\|_3 \|\vv_1 - \vv_2\|_2 \le \epsilon\|\vv_1 - \vv_2\|_6^2 + C\|\nabla \vv_1\|_3^2\|\vv_1 -\vv_2\|_2^2 .
\end{split}
\]
Then the application of the Sobolev embedding and Korn's inequality yields
\begin{equation} \label{uniQlast}
    {\rm Q} \le \epsilon\|\DD\vv_1 - \DD\vv_2\|_2^2 + C \|\nabla \vv_1\|_3^2 \|\vv_1 - \vv_2\|_2^2.
\end{equation}
Next, we shall write ${\rm P}$ as
\begin{equation} \label{uniN}
\begin{split}
    {\rm P} &= \int_\Omega \left( \SSS(c_1,\DD\vv_1)-\SSS(c_1,\DD\vv_2) \right):(\DD\vv_1 - \DD\vv_2)\dx \\
    & \hspace{4mm} + \int_{\Omega}\left(\SSS(c_1,\DD\vv_2) - \SSS(c_2,\DD\vv_2)\right) : (\DD\vv_1-\DD\vv_2) \dx \\
    &=: {\rm P_1} + {\rm P_2}.
\end{split}
\end{equation}
For the term ${\rm P_1}$, by the monotonicity assumption \eqref{P2}, we have
\begin{equation} \label{uniN1}
    {\rm P_1} \ge C \int_\Omega (1+|\DD\vv_1|^2+|\DD\vv_2|^2)^\frac{p(c_1)-2}{2} |\DD\vv_1 - \DD\vv_2|^2 \dx.
\end{equation}
Furthermore, by the mean value theorem, there exists $z$ between $c_1$ and $c_2$ such that
\[
\begin{split}
    {\rm P_2} 
    &= 2\int_\Omega p'(z) \nu(z,\DD\vv_2) \log(1+ |\DD\vv_2|^2) (c_1-c_2) \DD\vv_2 : (\DD\vv_1 - \DD\vv_2) \dx.
\end{split}
\]
Since $p \in W^{1,\infty}(\mathbb{R})$ and $\log(1+ |\DD\vv|^2) \le C (1+ |\DD\vv|^2)^\alpha$ for some arbitrarily small $\alpha>0$, we have
\begin{equation} \label{uniP2mid1}
\begin{split}
    |{\rm P_2}| &\le C \int_\Omega (1+ |\DD\vv_2|^2)^{\frac{p(z)-2}{2} + \frac{1}{2} + \alpha} |\DD\vv_1- \DD\vv_2||c_1-c_2| \dx \\
    &= C\int_\Omega (1+ |\DD\vv_2|^2)^{\frac{p(z)-2}{2} + \frac{1}{2} + \alpha - \frac{p(c_1)-2}{4}} (1+|\DD\vv_2|^2)^\frac{p(c_1)-2}{4}|\DD\vv_1 - \DD\vv_2||c_1 - c_2| \dx.
\end{split}
\end{equation}
We will estimate \eqref{uniP2mid1} for the cases of two and three-dimensional domains separately.

(2D case) By H\"older's inequality and Young's inequality, we see that
\[
    |{\rm P_2}| \le C \left( \int_\Omega (1+ |\DD\vv_2|^2)^{(2+\delta_1)(\frac{p^+}{2} - \frac{p^-}{4} + \alpha)} \dx \right)^\frac{2}{2+\delta_1} \|c_1 - c_2\|_\ell^2 + \epsilon \int_\Omega (1+ |\DD\vv_2|^2)^\frac{p(c_1)-2}{2} |\DD\vv_1 - \DD\vv_2|^2 \dx
\]
for some small $\delta_1>0$ and sufficiently large $\ell>2$ with $1= \frac{1}{2+\delta_1} + \frac{1}{\ell} + \frac{1}{2}$. From \eqref{2DgradetaL2} and the Sobolev embedding,
\[
\eta \in L^2(0,T;L^r(\Omega)), \quad 1 \le r < \infty,
\]
which implies
\[
\DD\vv_2 \in L^{p^-}(0,T; L^{\frac{p^-}{2} r} (\Omega)).
\]
Note further from \eqref{etainfty2D} that
\[
\DD\vv_2 \in L^\infty(0,T;L^{p^-}(\Omega)).
\]
If we choose large enough $r>0$, by the interpolation inequality, we can show that there exists small $\delta_2>0$ such that
\[
    \DD\vv_2 \in L^{2p^- - \delta_2}(Q_T).
\]
Therefore, if $p^+ < \frac{3}{2}p^-$, we have $(2+\delta_1)(p^+ - \frac{p^-}{2} + 2 \alpha )\le 2p^- - \delta_2$. Hence, by the Sobolev embedding, we obtain
\[
    |{\rm P_2}| \le F(t) \|\nabla (c_1-c_2)\|_2^2 + \epsilon\int_\Omega (1+ |\DD\vv_2|^2)^\frac{p(c_1)-2}{2} |\DD\vv_1 - \DD\vv_2|^2 \dx
\]
for some $F(t) \in L^1(0,T)$.\\

(3D case)
For the three-dimensional case, we have by H\"older's inequality and Young's inequality that
\[
    |{\rm P_2}| \le C\left( \int_\Omega (1+ |\DD\vv_2|^2)^{3(\frac{p^+}{2} - \frac{p^-}{4} + \alpha)} \dx \right)^\frac{2}{3} \|c_1- c_2\|_6^2 + \epsilon \int_\Omega (1+ |\DD\vv_2|^2)^\frac{p(c_1)-2}{2} |\DD\vv_1-\DD\vv_2|^2 \dx. 
\]
From \eqref{3DgradetaL2} together with the Sobolev embedding theorem, there holds
\[
    \eta \in L^2(0,T;L^6(\Omega)),
\]
which implies
\[
    \DD\vv_2 \in L^{p^-}(0,T;L^{3p^-}(\Omega)).
\]
Also by recalling \eqref{3DDVLinfty}, we have
\[
    \DD\vv_2 \in L^\infty(0,T;L^{p^-}(\Omega)).
\]
From the above estimates with the interpolation, we obtain
\[
    \DD\vv_2 \in L^{\frac{4}{3}p^-}(0,T;L^{2p^-}(\Omega)).
\]
Therefore, if $p^+ < \frac{7}{6}p^-$, we have $3(p^+ - \frac{p^-}{2} + 2\alpha) \le 2p^-$ which implies by the Sobolev embedding that
\[
    |{\rm P_2}| \le F(t) \|\nabla(c_1- c_2)\|_2^2 + \epsilon \int_\Omega (1+ |\DD\vv_2|^2)^\frac{p(c_1)-2}{2} |\DD\vv_1-\DD\vv_2|^2 \dx
\]
for some $F(t) \in L^1(0,T)$. All together, for both 2D and 3D cases, we have for some $F(t) \in L^1(0,T)$ that
\begin{equation} \label{uniN2}
    |{\rm P_2}| \le F(t) \|\nabla(c_1- c_2)\|_2^2 + \epsilon \int_\Omega (1+ |\DD\vv_2|^2)^\frac{p(c_1)-2}{2} |\DD\vv_1-\DD\vv_2|^2 \dx.
\end{equation}
Now we insert \eqref{uniN1} and \eqref{uniN2} into \eqref{uniN} to get
\begin{equation} \label{uniPlast3D}
    {\rm P} \ge C\int_\Omega(1+ |\DD\vv_1|^2+|\DD\vv_2|^2)^\frac{p(c_1)-2}{2} |\DD\vv_1 - \DD\vv_2|^2 \dx - F(t) \|\nabla (c_1-c_2)\|_2^2. 
\end{equation}
Then by combining \eqref{uniPQ}, \eqref{uniQlast} and \eqref{uniPlast3D}, we have
\begin{equation} \label{univelo3d}
    \frac{\rm d}{\dt} \|\vv_1 - \vv_2\|_2^2 + C \int_\Omega |\DD\vv_1 - \DD\vv_2|^2 \dx \le C \|\nabla \vv_1\|_3^2 \|\vv_1 - \vv_2\|_2^2 + F(t) \|\nabla (c_1-c_2)\|_2^2.
\end{equation}
Finally, from \eqref{uniconcentration} and \eqref{univelo3d}, we obtain
\begin{equation} \label{univandc}
\begin{split}
    &\frac{\rm d}{\dt} \|\nabla (c_1-c_2)\|_2^2 +  \frac{\rm d}{\dt} \|\vv_1 - \vv_2\|_2^2 + C\int_\Omega |\Delta (c_1-c_2)|^2 \dx + C \int_\Omega |\DD\vv_1 - \DD\vv_2|^2 \dx \\
    &\hspace{4mm}\le C \left(\|\vv_2\|_\infty^2 + F(t)\right) \|\nabla (c_1-c_2)\|_2^2  + C \left(\|\nabla c_1\|_\infty^2 + \|\nabla \vv_1\|_3^2\right) \|\vv_1 - \vv_2\|_2^2.
\end{split}
\end{equation}
Since $\vv_1, \vv_2 \in L^2(0,T;W^{2,2}(\Omega))$, by the Sobolev embedding,  we know that $\|\vv_2\|_\infty^2 \in L^1(0,T)$ and $\|\nabla \vv_1\|_3^2 \in L^1(0,T)$ for both 2D and 3D cases. Furthermore, based on the result \eqref{ct} and \eqref{ct3d}, we may conclude from Lemma \ref{cz} and the Sobolev inequality that $c_1 \in L^2(0,T; W^{1,\infty}(\Omega))$, which implies $\|\nabla c_1\|_\infty^2 \in L^1(0,T)$. Hence we can apply Gr\"onwall's inequality to \eqref{univandc} to get
\begin{equation} \label{uni3d}
    \nabla (c_1 - c_2) = 0 \quad \text{and} \quad \vv_1-\vv_2 = 0.
\end{equation}
Furthermore, by using the Poincar\'e inequality, we can conclude
\[
c_1 -c_2 = 0,
\]
which completes the proof.

\section{Conclusion}
In this paper, we have proved the existence of the unique strong solution for the system composed of the generalized non-Newtonian fluids with a concentration-dependent power-law index and the convection-diffusion equation in the two and three dimensions under the spatial periodicity condition. The main difficulty of the proof comes from the terms $\nabla c$ and $\pa_t c$ which appear when we differentiate the stress tensor in the Galerkin system. We first obtained the higher integrability of $\nabla c$ by testing the appropriate test function where the restriction $p^- \ge \frac{d+2}{2}$ is necessarily needed. To control $\pa_t c$, we have used the maximal regularity theory. Consequently, we obtained the existence of a strong solution with suitable regularity properties. Additionally, if we further assume the condition $p^+ < \frac{3}{2}p^-$ for $d=2$ and $p^+ < \frac{7}{6} p^-$ for $d=3$, we have also obtained the uniqueness of the strong solution.

An interesting future research direction is to obtain better regularity for the solution and extend this result to the degenerate case. In addition, from the computational point of view, performing the convergence analysis of the finite element approximation for the strong solution would also be intriguing, which will be addressed in the forthcoming paper.

\section*{Acknowledgements}
Kyueon Choi is supported by the National Research Foundation of
Korea Grant funded by the Korean Government (RS-2023-00212227).
Kyungkeun Kang is supported by the National Research Foundation of
Korea Grant funded by the Korean Government (RS-2024-00336346).
Seungchan Ko is supported by the National Research Foundation of
Korea Grant funded by the Korean Government (RS-2023-00212227).

\bibliography{references}
\bibliographystyle{abbrv}


\end{document}